\documentclass[10pt,a4paper]{article}
\usepackage[a4paper,top=2.5cm,bottom=2.5cm,left=1.5cm,right=1.5cm]{geometry}

\usepackage{amsmath}
\usepackage{color}
\usepackage{array}
\usepackage{ulem}
\usepackage{tabularx}
\usepackage{cancel}
\usepackage{subcaption}
\usepackage{tikz,pgfplots}
\usetikzlibrary{positioning,hobby,decorations.markings,arrows.meta,calc,fit}
\usepackage[numbers]{natbib}
\usepackage{url}
\usepackage{authblk}

\newcommand{\micr}{\,$\mu$m\,}

\newcommand{\adj}{\text{adj}}

\DeclareMathOperator{\e}{e}

\newcommand{\dod}{\!\cdot\!}

\newcommand{\sdot}{\!\cdot\!}

\newcommand{\nnum}{\nonumber \\}

\newcommand{\Gphi}{\nabla\vphi}

\newcommand{\Gtilphi}{\nabla\tilvphi}

\newcommand{\diver}{\text{div}}
\newcommand{\grad}{\text{grad}}

\newcommand{\tsym}{\text{sym}}

\newcommand{\bfens}{\bar{\nabla}}

\newcommand{\veps}{\varepsilon}
\newcommand{\vphi}{\varphi}

\newcommand{\cA}{{\cal A}}
\newcommand{\cC}{{\cal C}}

\newcommand{\cP}{{\cal P}}

\newcommand{\tDC}{\text{DC}}
\newcommand{\tAC}{\text{AC}}

\newcommand{\tA}{{\text A}}
\newcommand{\tD}{{\text D}}

\newcommand{\tT}{{\text T}}

\newcommand{\tL}{\text{L}}
\newcommand{\tK}{{\text K}}

\newcommand{\tU}{\text{U}}

\newcommand{\pk}{^{\text{PK}}}

\newcommand{\iA}{\int_{\text{A}}}
\newcommand{\iAz}{\int_{\text{A}_0}}

\newcommand{\dS}{\text{dS}}

\newcommand{\dOmega}{{\text d}\Omega}

\newcommand{\bfm}[1]{\mathbf{#1}}
\newcommand{\bsn}[1]{\boldsymbol{#1}}
\newcommand{\bfa}{\bfm{a}}
\newcommand{\bfb}{\bfm{b}}
\newcommand{\bfc}{\bfm{c}}
\newcommand{\bfe}{\bfm{e}}
\newcommand{\bff}{\bfm{f}}
\newcommand{\bfg}{\bfm{g}}

\newcommand{\bfn}{\bfm{n}}
\newcommand{\bfu}{\bfm{u}}

\newcommand{\bfs}{\bfm{s}}
\newcommand{\bfy}{\bfm{y}}

\newcommand{\bfx}{\bfm{x}}
\newcommand{\bfz}{\bfm{z}}
\newcommand{\bfone}{\bfm{1}}
\newcommand{\bfzero}{\bfm{0}}
\newcommand{\alphavec}{\boldsymbol \alpha}
\newcommand{\bfM}{\bfm{M}}
\newcommand{\bfK}{\bfm{K}}

\newcommand{\bfD}{\bfm{D}}
\newcommand{\bfA}{\bfm{A}}
\newcommand{\bfB}{\bfm{B}}
\newcommand{\bfC}{\bfm{C}}
\newcommand{\bfE}{\bfm{E}}
\newcommand{\bfF}{\bfm{F}}

\newcommand{\bfX}{\bfm{X}}
\newcommand{\bfY}{\bfm{Y}}

\newcommand{\bfR}{\bfm{R}}

\newcommand{\bsnR}{\bsn{R}}
\newcommand{\bsnD}{\bsn{D}}
\newcommand{\bsnC}{\bsn{C}}
\newcommand{\bsnG}{\bsn{G}}
\newcommand{\bsnH}{\bsn{H}}
\newcommand{\bsnQ}{\bsn{Q}}
\newcommand{\bsnT}{\bsn{T}}
\newcommand{\bfS}{\bfm{S}}
\newcommand{\bfV}{\bfm{V}}
\newcommand{\bfU}{\bfm{U}}
\newcommand{\bfW}{\bfm{W}}

\newcommand{\bfpsi}{\bsn{\psi}}
\newcommand{\bfsig}{\bsn{\sigma}}

\newcommand{\bfPhi}{\bsn{\Phi}}
\newcommand{\bfPsi}{\bsn{\Psi}}
\newcommand{\bfeps}{\bsn{\varepsilon}}

\newcommand{\tilu}{\tilde{u}}
\newcommand{\tils}{\tilde{s}}

\newcommand{\tilbfu}{\bfm{\tilde{u}}}

\newcommand{\nablaT}{\nabla^\tT}

\newcommand{\tilpsi}{\tilde{\psi}}
\newcommand{\tilvphi}{\tilde{\vphi}}

\newcommand{\da}{{\text d}a}

\newcommand{\ds}{{\text d}s}

\newcommand{\der}[2]{\frac{\partial #1}{\partial #2}}

\newcommand{\iu}{\pmb{\mathrm{i}}}

\newcommand{\tgreen}[1]{\textcolor{green!80!black}{#1}} 
\newcommand{\opcty}{1}

\definecolor{air_color}{RGB}{128, 223, 255}
\definecolor{potato_color}{RGB}{255, 219, 77}

\begin{document}

\title{Reduced order modelling of fully coupled electro-mechanical systems through invariant manifolds with applications to microstructures}

\author[1]{Attilio Frangi}
\author[1]{Alessio Colombo}
\author[2]{Alessandra Vizzaccaro}
\author[3]{Cyril Touzé}

\affil[1]{\small{Department of Civil and Environmental Engineering, Politecnico di Milano, P.za Leonardo da Vinci 32, 20133 Milano, Italy}}

\affil[2]{College of Engineering, Mathematics and Physical Sciences, University of Exeter, Exeter, United Kingdom}

\affil[3]{Institute of Mechanical Sciences and Industrial Applications (IMSIA), ENSTA Paris - CNRS - EDF - CEA, Institut Polytechnique de Paris, 828 Boulevard des Mar\'echaux, 91761 Palaiseau, France}

\date{\vspace{-5ex}}

\maketitle

\begin{abstract}
This paper presents the first application of the direct parametrisation method for invariant manifolds 
to a fully coupled multiphysics problem involving the nonlinear vibrations of deformable structures
subjected to an electrostatic field.
The formulation proposed is intended for model order reduction of electrostatically actuated resonating Micro-Electro-Mechanical Systems (MEMS). 
The continuous problem is first rewritten in a manner that can be directly handled by the parametrisation method, which relies upon automated asymptotic expansions. A new mixed fully Lagrangian formulation is thus proposed which contains only explicit polynomial nonlinearities, which is then discretised in the framework of finite element procedures.
Validation is performed on the classical parallel plate configuration, where different formulations using either the general framework, or an approximation of the electrostatic field due to the geometric configuration selected, are compared. Reduced-order models along these formulations are also compared to full-order simulations operated with a time integration approach.
Numerical results show a remarkable performance both in terms of accuracy and wealth of nonlinear effects that can be accounted for. In particular, the transition from hardening to softening behaviour of the primary resonance while increasing the constant voltage component of the electric actuation, is recovered. Secondary resonances leading to superharmonic and parametric resonances are also investigated with the reduced-order model.
\end{abstract}


\section{Introduction}

Micro-electro-mechanical systems (MEMS) have
increasingly attracted considerable interest due to their small size,
high reliability, and low power consumption. 
MEMS accelerometers, gyroscopes, pressure sensors, micromirrors, magnetometers, microphones and many others are now essential components in many devices of everyday life~\cite{corigliano2018mechanics}.
The analysis of the dynamical behaviour of MEMS devices requires time-dependent, nonlinear, multiphysics models including electromagnetics, piezoelectricity, and fluid-structure interaction. Microdevice engineering faces intricate geometries and is burdened by uncertainties regarding material parameters and fabrication imperfections. This makes traditional Full-Order simulation strategies extremely expensive, if not infeasible.

In MEMS resonators, for instance, frequency is highly sensitive to the level of the oscillation.
Indeed, while large actuation voltages might be desirable in a resonator
in order to decrease the phase noise and facilitate the design of the electronic control circuit, on the other hand, they can excite non-linear phenomena, 
both in the mechanical structures and in the electrostatic readout. 
In particular, a non-linear frequency response causes a shift in the resonance 
frequency of the oscillator~\cite{Nayfeh79,zega2020} and a consequent loss in terms of frequency 
stability of the device. 
For this reason, it is fundamental to be able to model the dynamic behaviour of such devices 
by taking into account all the possible sources of non-linearities during the design process.

Apart from real-time clocks, 
the literature on MEMS is rich in applications associated with nonlinear effects,
as recently reviewed in \cite{younis20}.
For instance, some MEMS like electrostatic micromirrors are intrinsically 
nonlinear~\cite{tie17} due to the large rotations of the reflecting surfaces.
In \cite{Baguet2019,7416140, GRENAT2022103903}, the behaviour of resonator arrays for mass 
sensing applications and of rate-integrating gyroscopes is addressed.
Moreover, different sources of non-linearities 
can be combined to achieve an overall linear response of the sensor, 
to compensate for temperature effects, instabilities or for sensing applications~\cite{8346730, elata2, Juillard_milos}. 
Non-linear phenomena like bistability~\cite{Ghayesh2017}, 
internal resonances~\cite{Houri2019, shoshani1, shoshani2},
self-induced parametric amplification~\cite{thomas13-APL,POLUNIN2017300,elata1} and frequency combs
\cite{seshia17,seshia18} have also been investigated in recent years, both theoretically and experimentally.

One can identify multiple sources of nonlinearities when transformations are no longer infinitesimal: geometric effects refer to the nonlinear evolution and coupling of stress components within a structure; electrostatic forces depend in an intrinsically nonlinear manner on geometrical gaps, just like gas dissipation.  
The combination of geometrical and electrostatic nonlinearities has attracted many investigations like the semianalytical approach discussed by \cite{younis20,younis03,younis05,younisbook,JUILLARD2015}.
The application of Von K{\'a}rm{\'a}n beam theory and analytical expressions of the electrostatic forces to parallel plates capacitors allows one to simulate the classical transition from hardening to softening behaviour
in the frequency response function, together with superharmonic and sub-harmonic resonances. However, the extension of these approaches to more general configurations is problematic.

In other simulations of coupled electromechanical problems \cite{aluru02,aluru04,batra06,batra07} the electrostatic field was computed resorting to the Boundary Element approach, with a heavy computational burden.

In recent years, several coupled electromechanical problems have been addressed in the literature using a reduction technique based on implicit condensation. This method, sometimes also denoted as the {\it applied force method}, relies on applying prescribed static forces to the structure to infer the nonlinear stiffness of a Reduced Order Model (ROM)~\cite{Hollkamp2008}. It has been first applied to MEMS in \cite{FRANGI2019}
and then extended in \cite{gobat2021reduced} to simulate 
a MEMS gyroscope test-structure exhibiting 1:2 internal resonance,
in \cite{zega2020} to present refined modelling of non-linearities in MEMS resonators, and in \cite{givois21-CS} to deal with coupled piezoelectric structures. However, this approach is known to have limitations that are intrinsically linked to the static underlying assumptions used in its derivation, see {\it e.g.} \cite{YichangICE,HallerSF,NicolaidouIceKE}. On the other hand, as explained for example in \cite{ReviewROMGEOMNL}, nonlinear reduction techniques based on the invariant manifold theory don't have such limitations and thus should be preferred for dynamical problems.

Model order reduction techniques using nonlinear normal modes (NNM) defined as invariant manifolds of the phase space have now a long history since the first works on the subject~\cite{ShawPierre91,ShawPierre93,touze03-NNM,TOUZE:JSV:2006}. In recent years, numerous developments have generalized the method to arbitrary order expansions that can be directly applied to finite element models, using the parametrisation method of invariant manifold~\cite{Cabre3,Haro}. For the case of damped systems, the uniqueness of the invariant manifold (NNM) that is tangent to its linear counterpart at the origin has been demonstrated and named spectral submanifold (SSM)~\cite{Haller2016}. Using a normal form style in the parametrisation method, direct computation of SSM-based ROMs has been shown in \cite{PONSIOEN2020,JAIN2021How,li2021periodic} and publicly released in the code \texttt{SSMtools}~\cite{ssmtools}. Direct computation using the normal form approach for structural mechanics has been developed in \cite{artDNF2020,AndreaROM,APL23} with special emphasis on using the method in a non-intrusive manner and tackling internal resonances. A direct application of the parametrisation method for systems featuring geometric nonlinearity, discretised by the FE procedure, has also been developed in \cite{vizza21high,opreni22high}, where applications to MEMS structures (microbeams, micromirrors) clearly underlined the potentiality of the technique, which drastically reduces the computational burden without losing any accuracy thanks to arbitrary order developments. An original treatment of the non-autonomous term is also proposed in \cite{vizza2023superharm}, that generalizes previous first-order developments used in \cite{opreni22high,JAIN2021How}, and allows dealing with {\it e.g.} superharmonic resonance. All these developments are embedded in the open-source code \texttt{MORFE} (Model order reduction for finite element problems)~\cite{morfe}.
The direct parametrisation method for invariant manifold (DPIM) has shown its great versatility and has been used to compute efficient ROMS also for wing-like structures~\cite{JAIN2021How},  viscoelastic beams with gyroscopic force~\cite{li2021periodic}, or rotating structures with centrifugal effect~\cite{Martin:rotation}. For MEMS problems with weakly coupled physics, it has also been used to take into account piezoelectric actuation~\cite{opreniPiezo} for devices where the polarization history is imposed. 

The present article aims to extend this reduction technique to 
the first truly multiphysics application, namely fully coupled electromechanical problems, and to electrostatically actuated MEMS structures in particular, which still represent the vast majority of such devices. While the DPIM works perfectly well with geometrical nonlinearities that can be easily expressed under a polynomial form, the main problem of the electromechanical problem is that it introduces non-polynomial nonlinearities in the formulation,
see {\it e.g.} \cite{aluru02,aluru04,batra07}, where
a Lagrangian formulation on the initial configuration was proposed. 
Since our goal here is to propose an efficient reduction technique based on the DPIM, the equations of motion are first reworked in order to express the semi-discretised finite element (FE) problem as differential-algebraic equations (DAE) containing only polynomial nonlinearities. Such a formulation is then compatible with the general derivation of the DPIM for non-autonomous problems discussed in \cite{vizza2023superharm}.
Consequently, an efficient and accurate ROM strategy can be derived for arbitrary order solutions of forced electromechanical systems, with direct application to MEMS design.

The paper is organized as follows. Section~\ref{sec:pbformulation} presents the general formalism for a generic electromechanical problem involving a deformable body and rigid electrodes. Particular attention is devoted to rewriting such a problem in the original configuration and introducing an auxiliary field in order to express the nonlinearities under polynomial form. Section~\ref{sec:solution} discusses the FE discretisation and summarizes the semi-discrete problem to solve. The reduction technique based on the DPIM is then presented starting from the semi-discrete DAE in Section~\ref{sec:dpim}. The configuration where the structure and the electrode are two parallel plates is then addressed in Section~\ref{sec:PPF}, and the resulting simplifications that can be used to express the electrostatic forces in such case are recalled. Finally, numerical results for a two-dimensional parallel plate problem are collected in Section~\ref{sec:results}.  Pull-in verifications, primary resonance solution with hardening/softening behaviour, and secondary resonances are investigated and comparisons between full-order simulations and ROMs are reported, showing the excellent accuracy provided by the reduction technique.


\section{Problem formulation}
\label{sec:pbformulation}

In this section, the equations of motion for the coupled nonlinear electromechanical problem are derived, and a particular emphasis is put on rewriting them in such a way that they can be readily solved by the direct parametrisation method for invariant manifold. This leads to adding new unknown fields in order to express the semi-discrete problem resulting from the finite element discretisation as a differential-algebraic equation (DAE)
with only polynomial nonlinearities.

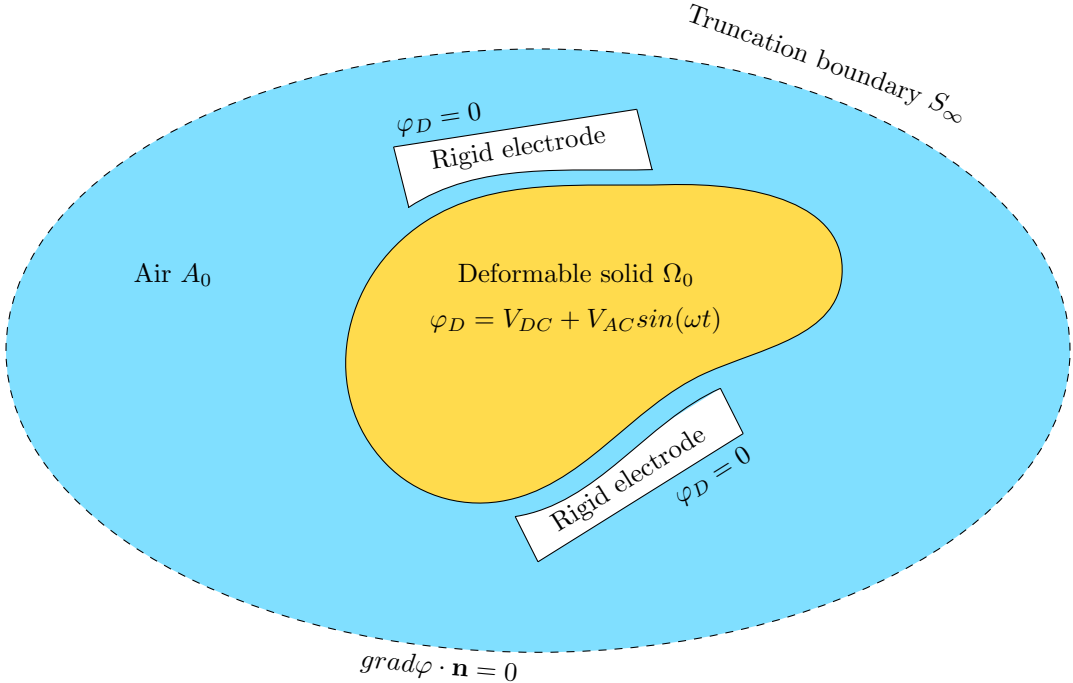
\begin{figure}[ht]
\centering
\begin{tikzpicture}[use Hobby shortcut]
    
    \node at (-0.5,-2) (origin) {};
    
    \node at (-1,-4) (a) {A};
    \node at (1.7,-2.3)  (b) {B};
    \node at (3.5,-1)   (c) {C};
    \node at (1.2,0.2)   (d) {D};
    \node at (-2,-0.3) (e) {E};
    \node at (-2.5,-3.5) (f) {F};

    \coordinate (shift1) at (0.2,-0.2);
    \node (a_l) at ($(a) + (shift1)$) {};
    \node (b_l) at ($(b) + (shift1)$) {};
    \node (c_l) at ($(c) + (shift1)$) {};
    \node (d_l) at ($(d) + (shift1)$) {};
    \node (e_l) at ($(e) + (shift1)$) {};
    \node (f_l) at ($(f) + (shift1)$) {};

    \coordinate (shift2) at (0.3,-0.6);
    \node (a_ll) at ($(a_l) + (shift2)$) {};
    \node (b_ll) at ($(b_l) + (shift2)$) {};

    \coordinate (shift3) at (-0.2,0.2);
    \node (a_u) at ($(a) + (shift3)$) {};
    \node (b_u) at ($(b) + (shift3)$) {};
    \node (c_u) at ($(c) + (shift3)$) {};
    \node (d_u) at ($(d) + (shift3)$) {};
    \node (e_u) at ($(e) + (shift3)$) {};
    \node (f_u) at ($(f) + (shift3)$) {};

    \coordinate (shift4) at (-0.2,0.8);
    \node (d_uu) at ($(d_u) + (shift4)$) {};
    \node (e_uu) at ($(e_u) + (shift4)$) {};

    \filldraw[dashed,color=air_color,draw=black,opacity=\opcty] (origin) ellipse (7 and 4);

    \fill[fill=white] (a_l.center) -- (b_l.center) -- (b_ll.center) -- (a_ll.center) -- cycle;
    \draw (b_l.center) -- (b_ll.center);
    \draw (b_ll.center) -- (a_ll.center);
    \draw (a_ll.center) -- (a_l.center);
    \begin{scope}
     \clip (a_l) rectangle (b_l);
        \filldraw[closed,color=air_color,draw=black,opacity=\opcty] (a_l) .. (b_l) .. (c_l) .. (d_l) .. (e_l) .. (f_l);
    \end{scope}

    \fill[fill=white] (d_u.center) -- (e_u.center) -- (e_uu.center) -- (d_uu.center) -- cycle;
    \draw (e_u.center) -- (e_uu.center);
    \draw (e_uu.center) -- (d_uu.center);
    \draw (d_uu.center) -- (d_u.center);
    \begin{scope}
     \clip (d_u) rectangle (e_u);
        \filldraw[closed,color=air_color,draw=black,opacity=\opcty] (a_u) .. (b_u) .. (c_u) .. (d_u) .. (e_u) .. (f_u);
    \end{scope}
    
    \filldraw[closed,color=potato_color,draw=black,opacity=\opcty] (a) .. (b) .. (c) .. (d) .. (e) .. (f);

    \pgfmathanglebetweenpoints{\pgfpointanchor{a_ll}{center}}{\pgfpointanchor{b_ll}{center}}
    \edef\myangle{\pgfmathresult}
    \node[above, sloped, rotate=\myangle] at ($(a_ll)!0.5!(b_ll)$) {Rigid electrode};
    \node[below, sloped, rotate=\myangle] at ($(a_ll)!0.8!(b_ll)$) {$\varphi_D = 0$};

    \pgfmathanglebetweenpoints{\pgfpointanchor{e_uu}{center}}{\pgfpointanchor{d_uu}{center}}
    \edef\myangle{\pgfmathresult}
    \node[below, sloped, rotate=\myangle] at ($(d_uu)!0.5!(e_uu)$) {Rigid electrode};
    \node[above, sloped, rotate=\myangle] at ($(d_uu)!0.8!(e_uu)$) {$\varphi_D = 0$};

    \node at ($(origin)+(0.5,1)$) {Deformable solid $\Omega_0$};
    \node at ($(origin)+(0.5,0.4)$) {$\varphi_D = V_{DC} + V_{AC} sin(\omega t)$};
    \node at ($(origin)+(-4.8,1)$) {Air $A_0$};
    \node[rotate=-20] at ($(origin)+(3.8,3.8)$) {Truncation boundary $S_\infty$};
    \node[rotate=-5] at ($(origin)+(-1.3,-4.2)$) {$grad \varphi \cdot \mathbf{n} = 0$};

\end{tikzpicture}
\caption{Example of the layout considered for the electromechanical problem, shown in the initial configuration.}
\label{fig:lay0}
\end{figure}

Let us consider a deformable solid, of initial configuration $\Omega_0$ 
as shown in Fig.~\ref{fig:lay0}.
The body is subjected to an electric field due to the interaction with a
set of rigid electrodes $\Omega_\tK$ 
and is otherwise immersed in the infinite empty space $A_0$.
The solid is considered a perfect conductor, a condition which is very well respected in silicon microsystems, with known given potential 
$\vphi^{\tD}=V_{\tDC}+V_{\tAC}\sin\omega t$, with $V_{\tAC}\ll V_{\tDC}$.
A uniform voltage $\vphi^{\tD}_\tK$ is imposed on each $\Omega_\tK$ (equipotential body)
and it is here assumed, for simplicity, that $\vphi^{\tD}_\tK=0$, 
although more general configurations can be easily considered.

To simplify the notation, $\vphi^{\tD}$ denotes the given voltage
on $\partial\Omega_0\cup S_0$, where $S_0$ is 
the collection of all the (fixed) surfaces of the electrodes. 


The problem is classically formulated in terms of a potential function 
$\varphi$ such that the electric field is $-\grad\varphi$,
and $\varphi$ is governed by the following set of
equations written in the current configuration:
\begin{alignat}{3}
\label{eq:strong1}
\diver\,\grad\varphi & =0, & \quad  & \text{in} \quad \, \tA,\\
\label{eq:strong2}
\varphi & =\vphi^{\tD}, & \quad & \text{on} \quad S_0, \\
\label{eq:strong3}
\grad\varphi\cdot\hat{\bfn} & =0, & \quad & \text{on} \quad S_{\infty},
\end{alignat}
where $\diver$ and $\grad$ denote operators with respect to the current coordinates
$\bfy$, $\tA$ is the current air domain, $S_{\infty}$ is a sufficiently distant truncation
boundary and the vector $\hat{\bfn}$ is the unit normal vector oriented 
from the solid into the air.
On the surface of a conductor, $\grad\vphi$ is orthogonal to the surface.
Alternatively, the problem can be formulated in a weak manner
as follows. 
\vspace*{.1cm}

Find $\vphi\in\cC_\varphi(\vphi^{\tD})$ such that:
\begin{align}
\label{eq:wf}
\iA \grad\tilvphi\sdot\grad\vphi\,\dOmega = 0,  \qquad \forall \tilvphi\in \cC_\varphi(0),
\end{align}
where $\cC_\varphi(\vphi^{\tD}$) denotes the space of sufficiently continuous functions 
respecting the Dirichlet boundary conditions given by Eq.~\eqref{eq:strong2}.
The solution to Eq.~\eqref{eq:wf} requires to mesh the air domain and to truncate it at $S_{\infty}$. 
Moreover, it is worth stressing that the domain $\tA$ is evolving according to the deformation of the solid.

The electrostatic pressure (force per unit surface),
exerted on the surface of the deformable body at an arbitrary point $\bfy$  reads:
\begin{align}
\label{eq:force}
\hat{\bfb} = \frac{\hat{\sigma}^2}{2\veps_0} \hat{\bfn}, 
\qquad
\hat{\sigma} = -\veps_0 \grad\vphi\cdot\hat{\bfn},
\end{align}
where $\hat{\sigma}$ is the surface charge on a conductor (Coulomb theorem).

%
%

The electrostatic forces induce vibrations in the structure
and the dynamical equilibrium is governed by the Principle of Virtual Power (PVP) enforced in the current unknown configuration:
\begin{equation}
\label{eq:PPVcurrent}
\int_{\Omega} \rho\ddot{\bfu}\cdot\tilbfu\,\dOmega + 
\int_{\Omega} \bfsig:\grad\tilbfu\,\dOmega = 
\int_{\partial\Omega} \hat{\bfb}\cdot\tilbfu\,\dS,
\qquad \forall\,\tilbfu\in\cC_u(\bfzero),
\end{equation}
where $\bfsig$ is the Cauchy stress tensor 
and $\hat{\bfb}$ are surface forces in the current configuration.
The test field $\tilbfu$  is defined over the space $\cC_u(\mathbf{0})$ of admissible functions, i.e. functions that are vanishing where Dirichlet boundary conditions are prescribed.
The first integral in Eq.\eqref{eq:PPVcurrent} expresses the power $\cP_a$ of inertia forces, while the second term is the opposite of the power of internal forces $\cP_i$,
and the right-hand side is the power of external forces $\cP_e$, such that the PVP can also be simply rewritten as:
\begin{equation}
\label{eq:PPVcurrent2}
\cP_a-\cP_i=\cP_e,
\qquad \forall\,\tilbfu\in\cC_u(\bfzero).
\end{equation}
It is important to remark that with this continuous weak formulation given by Eqs.~\eqref{eq:wf} and \eqref{eq:PPVcurrent}, the nonlinearities of the problem, stemming from both the mechanics (geometric nonlinearity) and the electric coupling, do not explicitly appear, being hidden in the fact that the current configuration is unknown. A classical method to make them explicit is to rewrite the problem in the original configuration.

\subsection{Formulation in the original configuration}
\label{sec:elecinitial}

The problem can be rewritten in the original configuration, using classical developments in continuum mechanics. For the sake of brevity, only the main results are reported here and the interested reader can find details of this derivation in appendix~\ref{app:MMCdetail}. The PVP~\eqref{eq:PPVcurrent} in the original configuration reads:
%
\begin{equation}
\label{eq:PPVmech}
\int_{\Omega_0} \rho_0\ddot{\bfu}\cdot\tilbfu\,\dOmega_0 + 
\int_{\Omega_0} \bfsig\pk[\bfu]:(\bfeps[\tilbfu]+\bfens[\bfu,\tilbfu])\,\dOmega_0 
=\int_{\partial\Omega_0} \bfb\cdot\tilbfu\,\dS_0
\qquad \forall\,\tilbfu\in\cC_u(\bfzero),
\end{equation}
where $\bfeps$ and $\bfens$, defined in Appendix~\ref{app:MMCdetail}, stand for:
\[
\bfeps[\tilbfu] =\frac{1}{2}(\nabla\tilbfu+\nablaT\tilbfu) 
\qquad
\bfens[\bfu,\tilbfu] = \frac{1}{2}(\nablaT\bfu\cdot\nabla\tilbfu+\nablaT\tilbfu\cdot\nabla\bfu), 
\]
%


%
$\nabla$ denotes the gradient with respect to the initial coordinates $\bfx$,
and $\bfy=\bfx+\bfu$.
In Eq.~\eqref{eq:PPVmech}, 
$\bfsig\pk$ denotes the second Piola-Kirchhoff stress tensor. 
Since for the class of problems at hand strains are small even if transformations are non-infinitesimal, the Kirchhoff constitutive law holds, yielding:
%
\begin{equation}
\label{eq:const}
\bfsig\pk =\cA:\bfe = \bfe:\cA,
\end{equation}
where $\cA$ is the tensor of elastic constants and $\bfe$ is
the Green-Lagrange strain tensor:
\begin{equation}
\label{eq:gl}
\bfe[\bfu]=\frac{1}{2}(\nabla\bfu+\nablaT\bfu+\nablaT\bfu\cdot\nabla\bfu).
\end{equation}
Let us now detail the derivation of the electrostatic force in the initial configuration. Since it is a follower force, the calculations are not straightforward and have been thus reported in Appendix~\ref{app:surface}. 
In particular, using Eq.~\eqref{eq:Nds}, one has:
\begin{align}
\label{eq:Pe}
\cP_e  & =\int_{\partial\Omega}\frac{\hat{\sigma}^2}{2\veps_0} (\tilde{\bfu}\cdot\hat{\bfn})\dS
\\
& =\frac{\veps_0}{2}\int_{\partial\Omega_0} \sigma^2
\left(\bfn +  \frac{1}{J_S}(\bfg_1\wedge\bfu_{,2}+\bfu_{,1}\wedge\bfg_2 + \bfu_{,1}\wedge\bfu_{,2})\right)\cdot\tilde{\bfu}\,\dS_0,
\nonumber
\end{align}
where $J_S$ is the surface jacobian, see Eq.\eqref{a2},
$\bfu_{,\alpha}$ denotes the derivative with respect to 
the surface coordinate $a_\alpha$ used to represent the FE surface elements,
$\bfg_\alpha$ are surface covariant vectors and 
$\sigma(\bfx)=\hat{\sigma}(\bfy(\bfx))$.
Eq.~\eqref{eq:Pe} indirectly provides an explicit expression for $\bfb$  in Eq.~\eqref{eq:PPVmech}:
\[
\bfb=\sigma^2
\left(\bfn +  \frac{1}{J_S}(\bfg_1\wedge\bfu_{,2}+\bfu_{,1}\wedge\bfg_2 + \bfu_{,1}\wedge\bfu_{,2})\right).
\]
In order to propose a unified formulation expressed in the initial configuration, the electrostatics equilibrium, Eq.~\eqref{eq:wf}, also needs to be rewritten, a step that is rarely done in the literature and is detailed next.
One of the problems to be solved is linked to the fact that the deformation of the air domain $\tA$ follows the displacement of the structure $\bfu$. Since this displacement is defined only for the structural nodes, the air domain deformation is strictly speaking only defined on $\partial\Omega$. One method to cope with this issue is to extend smoothly the structural displacement vector $\bfu$ to $\tA$ to represent the {\it deformation} of the air domain.
This extension, denoted $\bfs$, is indeed fictitious and not unique.
In a FEM sense, one could a priori set the nodal values to zero in the air domain 
immediately outside the conductors, but this leads to strong gradients, 
thus generating a poor solution quality. 
The option adopted in the present investigation is to solve a Dirichlet problem for each displacement component, i.e.\ find $s_i\in \cC_s(u_i)$, $\forall i \in [1,d]$:
\begin{align}
\label{eq:wref_laplace}
\iAz \nabla\tils\sdot\nabla s_i\,\dOmega_0 = 0   
\quad \forall \tils\in \cC_s(0),
\end{align}
where the space $\cC_s(u_i)$ collects the functions such that 
$s_i=u_i$ on $\partial\Omega_0$, $s_i=0$ on $S_0\cup S_\infty$ and $d$ is the dimension of the problem, $d=2$ if the problem is 2-dimensional, $d=3$ if it is 3-dimensional.
An example of this extension procedure is presented in Section \ref{sec:results},
where a clamped-clamped beam faces a rigid electrode.
Other choices that would be computationally more efficient will be commented on in Section~\ref{sec:remarksBL}.
Once $\bfs$ has been defined, the electrical equilibrium in the initial configuration reads
\begin{align}
\label{eq:wref}
\iAz \Gtilphi\sdot\bfc^{-1}[\bfs]\sdot\Gphi J\dOmega_0 = 0,   
\qquad \forall \tilvphi\in \cC_\varphi(0),
\end{align}
where  $J=\det\bff[\bfs]$ is the Jacobian of the transformation:
\begin{align}
\bff[\bfs] & =\bfone+\nabla\bfs,
\\
\bfc[\bfs] & =\bff^\tT[\bfs]\cdot\bff[\bfs]=\bfone+\nabla\bfs+\nabla^\tT\bfs+
\nabla^\tT\bfs\sdot\nabla\bfs.
\end{align}
%
The main concern with this formulation is that $\bfc^{-1}$ is not polynomial in $\bfs$. The next section aims to circumvent this difficulty by adopting a new formulation of the continuous problem.

\subsection{Mixed formulation with polynomial nonlinearities}
\label{sec:addeq}

The goal of this Section is to show how the introduction of an additional field allows the expression of the problem to be solved with polynomial nonlinearities only. 
To remove the non-polynomial terms
from Eq.~\eqref{eq:wref}, an additional vector field $\bfpsi$ can be introduced as:
\begin{equation}
\label{eq:new}
\bfpsi=\Gphi\cdot\bff^{-1}[\bfs]=\bff^{-\tT}[\bfs]\cdot\Gphi, \qquad
\Gphi=\bff^\tT[\bfs]\cdot\bfpsi=\bfpsi\cdot\bff[\bfs].
\end{equation}
Remark that $\bfpsi$ corresponds to the gradient in actual coordinates 
and is strictly related to the electric field.
The space of $\bfpsi$ might be for instance that of piecewise (discontinuous) linear vectors
(element by element).
We enforce Eq.\eqref{eq:new} in a weak form, using as test functions all possible
$\tilde{\bfpsi}\in \cC_{\psi}$:
\begin{align}
\label{eq:weak}
\iAz \tilde{\bfpsi}\sdot\bigl(\bff^\tT[\bfs]\cdot\bfpsi-\Gphi)\dOmega_0 = 0,   
\qquad \forall\, \tilde{\bfpsi} \in \cC_{\psi},
\end{align}
while Eq.~\eqref{eq:wref} now becomes:
\begin{align}
\label{eq:wref2}
\iAz \Gtilphi\cdot(J\bff^{-1}[\bfs])\cdot\bfpsi\,\dOmega_0 = 0.,  
\qquad \forall \tilvphi\in \cC_\varphi(0).
\end{align}
As $J\bff^{-1}$ is the adjoint of $\bff$, i.e.\ the transpose of the cofactor matrix of $\bff$, we will use the notation $\adj(\bff)$ in what follows.
Recall that the coefficient $(i,j)$ of the cofactor matrix of 
$\bff$ is $(-1)^{i+j}M_{ij}$, where $M_{ij}$ is the minor of the coefficient.
Hence, each coefficient of $\adj(\bff)[\bfs]$ is a quadratic polynomial in $\bfs$ in a 3D space. On the other hand, for 2D space problems, each coefficient of $\adj(\bff)[\bfs]$ is linear in $\bfs$. Such a simplification is helpful in solving planar problems and will be used in Section~\ref{sec:results}.

Since on the surface of the conductors $\grad\vphi$ is parallel to $\hat{\bfn}$, the following relationships hold:
\[
\sigma\hat{\bfn} = -\veps_0 \Gphi\cdot\bff^{-1},
\qquad
\sigma^2=\veps^2_0 \|\bfpsi\|^2,
\]
and the virtual power of the electrostatic pressure 
in the reference configuration can be written as a function of the additional vector field $\bfpsi$ in a formulation without non-polynomial nonlinearities as:
\begin{align}
\label{eq:powerEM}
\cP_e =
\frac{\veps_0}{2}\int_{\partial\Omega_0} \|\bfpsi\|^2
\left(\bfn+\frac{1}{J_S}(\bfg_1\wedge\bfu_{,2}+\bfu_{,1}\wedge\bfg_2 + \bfu_{,1}\wedge\bfu_{,2})\right)\sdot\tilde{\bfu}\,\dS_0.
\end{align}
%

\subsection{Summary of the General Formulation}
\label{sec:GF}

To summarize, we have proposed a formulation of the coupled electromechanical problem 
consisting of a system of four equations
which can be expressed as follows.
\vspace{.1cm}

Find $\bfu\in\cC_u(0)$, $\bfs\in\cC_s(\bfu)$,
$\varphi\in \cC_\varphi(\varphi^D)$,  $\bfpsi\in\cC_\psi$
such that:

\begin{itemize}
\item
the mechanical equilibrium:
\begin{multline}
\label{eq:GFmech}
\int_{\Omega_0} \rho_0\ddot{\bfu}\cdot\tilbfu\,\dOmega_0 + 
\int_{\Omega_0} \bfsig\pk[\bfu]:(\bfeps[\tilbfu] + \bfens[\bfu,\tilbfu])\,\dOmega_0 
\\
=\frac{\veps_0}{2}\int_{\partial\Omega_0} \|\bfpsi\|^2
\left(\bfn+\frac{1}{J_S}(\bfg_1\wedge\bfu_{,2}+\bfu_{,1}\wedge\bfg_2 + \bfu_{,1}\wedge\bfu_{,2})\right)\sdot\tilde{\bfu}\,\dS_0,
\quad \forall \tilbfu \in \cC_u(0),
\end{multline}
\item
the extension problem, $\forall i \in [1,d]$:
\begin{align}
\label{eq:GFext}
\iAz \nabla\tils\sdot\nabla s_i \,\dOmega_0 = 0,   
\quad \forall \tils\in \cC_s(u_i),
\end{align}
\item
the compatibility equation:
\begin{align}
\label{eq:GFcomp}
\iAz \tilde{\bfpsi}\sdot\bigl(\bfpsi+\nabla^\tT\bfs\cdot\bfpsi-\Gphi)\dOmega_0 = 0,   
\qquad \forall\, \tilde{\bfpsi} \in \cC_{\psi},
\end{align}
\item
the electrostatics equilibrium: 
\begin{align}
\label{eq:GFelet}
\iAz \Gtilphi\cdot(\bfpsi+\adj(\nabla\bfs)\cdot\bfpsi)\,\dOmega_0 = 0,   
\qquad \forall \tilvphi\in \cC_\varphi(0),
\end{align}
\end{itemize}
are all satisfied.
\vspace{.1cm}

The system of equations~\eqref{eq:GFmech}-\eqref{eq:GFelet} will be referred to as the Mixed General Formulation ({\sc mgf}) in the rest of the paper, as it mixes standard kinematic unknowns like displacements and new {\it force} variables like $\bfpsi$.

It provides a very general framework to solve for coupled electromechanical problems, which also fits the requirements of the ROM strategy using the DPIM. Indeed, the problem with non-polynomial nonlinearities has been replaced by a DAE containing only quadratic and cubic terms, at the price of introducing an auxiliary field, which results in an added compatibility equation. Section~\ref{sec:dpim} will show how the general reduction strategy provided by the DPIM can be applied to the semi-discrete problem resulting from a FE discretisation of these equations.


\subsubsection{Remark on computational efficiency}
\label{sec:remarksBL}

From the inspection of the {\sc mgf} it is clear
that the introduction of the two fields $\bfs$ and $\bfpsi$
is a considerable burden from the computational point of view. 
However, this can be drastically reduced
by defining, in the original configuration, a boundary layer $\tL_0$ of suitable thickness around the deformable structure
such that $\bfs$ decays smoothly to zero from $\partial\Omega_0$ 
to the outer boundary $\partial\tL_0$ of the layer.
In $\tA_0  \backslash \tL_0$, $\bfs$ vanishes, $\bfc[\bfs]=\bfone$, there is no
need to introduce $\bfpsi$ and the original trivial electrostatic formulation is recovered.

\subsubsection{Simplifications in 2D}
\label{sec:2D}

Since the results analysed in Section \ref{sec:results}
are obtained using a 2D implementation,
 the simplifications pertaining to this specific case are here detailed.
Solids are extruded along the $z$ direction,
the first surface coordinate $a_1$ coincides with $z$ and hence 
$\bfg_1=\bfe_z$ and $\bfu_{,z}=\bfm{0}$. 
The outer surfaces of the solids are represented in 2D by the lines
which are their intersection with the $z=0$ plane.
The virtual power $\cP_e$ of electrostatic forces hence becomes:
\begin{align}
\label{eq:PPVS2d}
\cP_e=\frac{\veps_0}{2}\int_{\partial\Omega_0} \|\bfpsi\|^2\left(\bfn +  
\frac{\bfe_z\wedge\bfu_{,a}}{J_S}\right)\sdot\tilbfu\,\ds,
\end{align}
where $a$ is an abscissa running along the line-boundary $\partial\Omega_0$,
$\bfu_{,a}$ is the directional derivative of $\bfu$ with respect to $a$
and $\ds=J_S\da$.
Moreover, $\adj(\nabla\bfs)$ is only linear in $\bfs$.

\section{Solution procedure}
\label{sec:solution}

In coupled electromechanical problems, vibrations develop around
an initial static solution associated with the presence of the
constant voltage $V_\tDC$. This is known to have important effects
on the nonlinear response of the structure, as will be evidenced in
Section~\ref{sec:results}.

In the first step, the static solution is computed 
using a Newton-Raphson solver  linearizing the governing equations
around the current guess of the solution and generating the so-called
{\it tangent stiffness matrix} at every iteration. It is worth stressing that 
this matrix, at convergence, will also represent the linear part of the 
expansion in the DPIM procedure, see Section~\ref{sec:dpim}, that, together with the
mass matrix, defines the eigenvalue problem. 

We will henceforth denote with $\bfu_0,\bfs_0,\bfpsi_0,\vphi_0$ the static solution
and with $\bfu_0+\bfu,\bfs_0+\bfs,\bfpsi_0+\bfpsi,\vphi_0+\vphi$ the 
total solution without changing notation, for the sake of simplicity.

In the next section, we will inject this decomposition into the governing equations
and partition all the terms into linear, quadratic and cubic contributions, in order to prepare 
the path for the discretisation and subsequent application of the DPIM to obtain a reduced-order model.
Constant terms only define the static solution and will be disregarded.

\subsection{FEM discretisation}

The discretised version of the equations is here provided according to standard FE procedures, by focusing on two-dimensional planar problems, for the sake of simplicity. Note that there is no difference between the 2D and 3D problems from the theoretical point of view, nevertheless, 
to obtain a rapid proof of concept of the reduction method, the presentation and the results shown in Section~\ref{sec:results} will only consider 2D problems.
All constant terms associated with the initial static solution are removed.
In the present implementation, $\bfu$ and $\bfs$ are interpolated with quadratic triangles,
as well as $\vphi$. On the contrary, $\bfpsi$ is modelled as piecewise-linear, discontinuous from triangle to triangle.

\paragraph{Mechanical equilibrium}

We analyse here the terms associated with the virtual power of electrostatic forces
in Eq.~\eqref{eq:GFmech}, while we simply provide the final result for 
the more classical mechanical terms, since they have already been considered
elsewhere, see {\it e.g.} \cite{opreniPiezo}, and are briefly discussed in Appendix \ref{app:MMCdetail}. 
If $\bfn_0$ denotes the normal to the surface deformed by the static solution $\bfu_0$, and $J_{S0}$ is the corresponding surface Jacobian, 
\begin{align}
\label{eq:PPVS2dS}
\cP_e =
& 
\underbrace{
\frac{\veps_0}{2}\int_{\partial\Omega_0} \tilde{\bfu}\sdot\left(2\bfn_0(\bfpsi_0\sdot\bfpsi) +\|\bfpsi_0\|^2 \frac{\bfe_z\wedge\bfu_{,a}}{J_{S0}}\right)\ds  
}_{\text{linear}}
\nnum
& + 
\underbrace{
\frac{\veps_0}{2}\int_{\partial\Omega_0} \tilde{\bfu}\sdot\left(\bfn_0\|\bfpsi\|^2 +
2\frac{\bfe_z\wedge\bfu_{,a}}{J_{S0}}(\bfpsi_0\sdot\bfpsi)\right)\ds
}_{\text{quadratic}}
+
\underbrace{
\frac{\veps_0}{2}\int_{\partial\Omega_0}
\tilde{\bfu}\cdot\frac{\bfe_z\wedge\bfu_{,a}}{J_{S0}}\|\bfpsi\|^2\ds,
}_{\text{cubic}}
\end{align}
which generates the discretised contributions:
\begin{equation}
\cP_e= \tilde{\bfU}^\tT\left(\bfR_u\bfU + \bfR_\psi\bfPsi+ 
\bsnR_{\psi\psi}(\bfPsi,\bfPsi)+\bsnR_{\psi u}(\bfPsi,\bfU)
+\bsnR_{\psi\psi u}(\bfPsi,\bfPsi,\bfU)\right).
\end{equation}

Here $\bfU$ and $\tilde{\bfU}$ are vectors collecting the unknown and test nodal values, respectively, for the $\bfu$ and $\tilbfu$ fields.
Similarly, $\bfPsi$ collects the nodal values of the $\bfpsi$ discretisation,
element by element. Other upright bold symbols, like $\bfK$, denote matrices,
while italic bold symbols, like $\bsnR$, denote quadratic or cubic nonlinear operators.
When added to the $\cP_a$ and $\cP_i$ terms, the final semi-discrete equation reads:
\begin{align}
\label{eq:discrM2}   
&\bfM \ddot{\bfU}+(\alpha \bfM + \beta \bfK)\dot{\bfU}+(\bfK-\bfR_u)\bfU-\bfR_\psi\bfPsi+ 
\bsnG(\bfU,\bfU)+ \bsnH(\bfU,\bfU,\bfU)+
\nnum
&-\bsnR_{\psi\psi}(\bfPsi,\bfPsi)-\bsnR_{\psi u}(\bfPsi,\bfU) - \bsnR_{\psi\psi u}(\bfPsi,\bfPsi,\bfU)=\bfzero,
\end{align}
where a mechanical dissipation term has been added in the form of Rayleigh damping.

\paragraph{Extension equation}

The extension equation is linear and is the easiest one to discretise. It is written for all the components
$s_i$ at once by collecting in $\bfS$  the unknown nodal values of the interpolation for 
$\bfs$:
\begin{align}
\label{eq:discrE2}
\bfE_s\bfS + \bfE_u\bfU = \bfzero,
\end{align}
where $\bfE_s$ is a square positive definite matrix and $\bfE_u\bfU$ accounts for 
the Dirichlet boundary conditions imposing that $\bfs=\bfu$ on $\partial\Omega_0$.

\paragraph{Compatibility equation}

The compatibility equation, setting $\bff_0=\bfone+\nabla\bfs_0$, becomes:
\begin{align}
\label{eq:GFcomp0}
\underbrace{
\iAz \tilde{\bfpsi}\sdot\bigl(\bff^\tT_0\cdot\bfpsi+\nabla^\tT\bfs\cdot\bfpsi_0-\Gphi)\dOmega_0
}_{\text{linear}}
+
\underbrace{
\iAz \tilde{\bfpsi}\sdot\nabla^T\bfs\cdot\bfpsi\,\dOmega_0
}_{\text{quadratic}} 
= 0,    
\end{align}
which generates the discretised form:
\begin{align}
\label{eq:discrC2}
\bfC_u \bfU+\bfC_s \bfS+\bfC_\psi \bfPsi+\bfC_\varphi\bfPhi+\bsnC_{\psi w}(\bfPsi,\bfS)
= \bfC_V V_{\tAC}\sin(\omega t),
\end{align}
where $\bfPhi$ is the collection of the unknown nodal values for the $\varphi$ field.
It is worth remarking that $\bfC_u \bfU$  stems from the Dirichlet boundary conditions on $\bfs$,
while $\bfC_V$ is due to the imposed voltage 
which fixes the nodal values of $\varphi$ on $\partial\Omega_0$. 
It is also important to stress that $V_\tDC$ appears
in the above equations only through the static solution which has been computed as a first step.

\paragraph{Electrostatics equilibrium}

The last equation of the problem concerns the electrostatics equilibrium and reads
\begin{align}
\label{eq:ES2d}
\underbrace{
\iAz\Gtilphi\cdot\left(\adj(\bff_0)\cdot\bfpsi+\adj(\nabla\bfs)\cdot\bfpsi_0\right)\dOmega_0
}_{\text{linear}}
+
\underbrace{
\iAz \Gtilphi\cdot\adj(\nabla\bfs)\cdot\bfpsi\dOmega_0 = 0,   
}_{\text{quadratic}}
\end{align}
which then generates the discretised form:
\begin{align}
\label{eq:discrD2}
\bfD_u\bfU+\bfD_s\bfS+\bfD_\psi\bfPsi+ \bsnD_{\psi s}(\bfPsi,\bfS)
 =\bfzero.
\end{align}
%

\section{Model order reduction}
\label{sec:dpim}

Now that the semi-discretised problem has been derived, the model-order reduction technique using the direct parametrisation method for invariant manifolds (DPIM) can be directly applied. This nonlinear reduction technique is a simulation-free method that operates directly from the physical space, and computes arbitrary order nonlinear mappings and reduced dynamics along the selected invariant manifold. To apply the DPIM algorithm to this problem, the first step consists of writing the system as  first-order DAE. To cope with the second-order time derivative in Eq.~\eqref{eq:discrM2}, an additional velocity field $\bfV$ is introduced via $\bfM \dot{\bfU} = \bfM \bfV$. Combining this  additional equation with Eq.~\eqref{eq:discrM2}, Eq.~\eqref{eq:discrE2}, Eq.~\eqref{eq:discrC2}, and Eq.~\eqref{eq:discrD2}, one obtains:
\begin{subequations}\begin{align}
    & 
    \bfM \dot{\bfV} = -(\alpha\bfM+\beta\bfK)\bfV+
    (-\bfK+\bfR_u)\bfU+\bfR_\psi\bfPsi 
    -\bsnG(\bfU,\bfU) - \bsnH(\bfU,\bfU,\bfU)+
    \nnum
    &\phantom{\bfM \dot{\bfV} =} +\bsnR_{\psi\psi}(\bfPsi,\bfPsi)+
    \bsnR_{\psi u}(\bfPsi,\bfU) + \bsnR_{\psi\psi u}(\bfPsi,\bfPsi,\bfU)
\\
    & \bfM \dot{\bfU} = \bfM \bfV
\\
    & \bfzero = -\bfE_s\bfS  -\bfE_u\bfU
\\
    & \bfzero = -\bfC_u \bfU
    -\bfC_s \bfS -\bfC_\psi \bfPsi -\bfC_\varphi\bfPhi-\bsnC_{\psi s}(\bfPsi,\bfS)
    +\bfC_V V_{\tAC}\sin(\omega t)
\\
    & \bfzero = -\bfD_u\bfU -\bfD_s\bfS -\bfD_\psi\bfPsi -\bsnD_{\psi s}(\bfPsi,\bfS),
\end{align}
\label{eq:original_DAEsystem}
\end{subequations}
where the derivative terms have been brought to the left-hand side and the algebraic terms on the right-hand side of the equations.

We can collect all the variables of the problem into a single variable $\bfy$ defined as:
\begin{equation}
    \bfy = 
    \begin{bmatrix}
        \bfV\\
        \bfU\\
        \bfS\\
        \bfPsi\\
        \bfPhi 
    \end{bmatrix}
\end{equation}

and one gets a compact form for the problem in the sole $\bfy$ variable with nonlinearities up to cubic order, which reads:
\begin{align}
    \label{eq:ODE_original}
    {\bfB}
        \dot{\bfy}
    =
    {\bfA}
    \bfy
    +
    \bsnQ(\bfy, \bfy)
    +
    \bsnT(\bfy, \bfy, \bfy)
    +
    \varepsilon (\bfC_+ \e^{+\iu \omega t}
    +
     \bfC_- \e^{-\iu \omega t}).
\end{align}

In particular, the linear terms of the system are collected by means of the matrices ${\bfA}$ and ${\bfB}$, defined as:
\begin{equation}
    {\bfB} \doteq     \begin{bmatrix}
        \bfM & \bfzero & \bfzero & \bfzero & \bfzero \\
        \bfzero & \bfM & \bfzero & \bfzero& \bfzero  \\
        \bfzero & \bfzero & \bfzero & \bfzero & \bfzero  \\
        \bfzero & \bfzero & \bfzero & \bfzero & \bfzero \\
        \bfzero & \bfzero & \bfzero & \bfzero & \bfzero 
    \end{bmatrix},
\qquad
    {\bfA} \doteq     \begin{bmatrix}
          - (\alpha\bfM+\beta\bfK)
        & - {\bfK} + {\bfR}_u
        & \bfzero 
        & {\bfR}_\psi
        & \bfzero 
        \\
        \bfM & \bfzero & \bfzero & \bfzero & \bfzero 
        \\
        \bfzero & -\bfE_u & -\bfE_s & \bfzero & \bfzero 
        \\
        \bfzero & -\bfC_u & -\bfC_s & -\bfC_\psi & -\bfC_\phi\\
        \bfzero & -\bfD_u & -\bfD_s & -\bfD_\psi & \bfzero
    \end{bmatrix}.
\end{equation}

The polynomial nonlinearities are collected into two vectors separating the quadratic forces $\bsnQ$ from the cubic ones $\bsnT$, which are defined as:
\begin{equation}
    \bsnQ(\bfy, \bfy) \doteq
    \begin{bmatrix}
    - {\bsnG}(\bfU,\bfU) 
    + {\bsnR}_{\psi\psi}(\bfPsi,\bfPsi)
    + {\bsnR}_{\psi u}(\bfPsi,\bfU)
    \\ \bfzero \\ \bfzero \\ -\bsnC_{\psi s}(\bfPsi, \bfS) \\ -\bsnD_{\psi s}(\bfPsi, \bfS)
    \end{bmatrix}, \quad \bsnT(\bfy, \bfy, \bfy) \doteq
    \begin{bmatrix}
    - \bsnH(\bfU,\bfU,\bfU) 
    + \bsnR_{\psi\psi u}(\bfPsi,\bfPsi,\bfU)
    \\ \bfzero \\ \bfzero \\ \bfzero \\ \bfzero
    \end{bmatrix}.
\end{equation}


Lastly, the non-autonomous term has been defined as:
\begin{equation}
        \varepsilon\bfC_\pm 
    \doteq
    \pm \dfrac{1}{2\;\iu}
    \begin{bmatrix}
        \bfzero\\ \bfzero\\ \bfzero \\ \bfC_V \\ \bfzero
    \end{bmatrix}
    V_{\tAC}, \quad \mathrm{such} \; \, \mathrm{that:} \quad 
    \varepsilon (\bfC_+ e^{+\iu \omega t}
    +
     \bfC_- e^{-\iu \omega t})
    =
    \begin{bmatrix}
        \bfzero\\ \bfzero\\ \bfzero \\ \bfC_V \\ \bfzero
    \end{bmatrix}
    V_{\tAC}
    \sin(\omega t).
\end{equation}

Note that the meaning of $\varepsilon$ coincides with the amplitude of the non-autonomous term. Consequently, $\bfC_\pm$ has to be unitary in some sense. To meet this constraint, $\bfC_\pm$ will be normalised with respect to one of the eigenvectors of the system. 

The compact form of the problem around the fixed point given by Eq.~\eqref{eq:ODE_original} is the starting point for the non-autonomous DPIM algorithm introduced in \cite{vizza2023superharm}. The procedure performs a parametrisation of an invariant manifold tangent to a selected master subspace of the system, linearised around the given fixed point. The first step of the DPIM is to compute the solution of the eigenproblem for a complex conjugate pair (or more) that will define the master subspace.

Let us define the triplet $\lbrace \lambda_k, \bfY_k, \bfX_k \rbrace$ to be one eigenvalue, right eigenvector, and conjugate left eigenvector of the system respectively. They fulfil the following properties:
\begin{subequations}
    \begin{align}
        & (\lambda_k \bfB - \bfA) \bfY_k = \bfzero
        \\
        & (\lambda_k \bfB^\tT - \bfA^\tT) \bfX_k = \bfzero.
    \end{align}
\end{subequations}

The normalisation condition of the eigenvectors reads:
\begin{equation}
    \bfX_k^\tT \bfB \bfY_l = b_k \delta_{kl},  \quad \forall k,l
\end{equation}

Therefore, the normalisation condition for the forcing reads:
\begin{equation}
    || \bfX_k^\tT \bfC_\pm ||_2 = b_k.
\end{equation}

The master subspace is typically selected to be the subspace spanned by the slowest modes of the system, i.e. the eigenvectors related to the eigenvalues with the smallest real part. In the non-autonomous context, we are interested in reducing the system around resonance. The forcing frequency $\omega$ will be close to a complex conjugate pair of eigenvalues of the system, thus exciting one of the modes: $\lambda_k \approx +\iu\omega $ and $\bar{\lambda}_k \approx -\iu\omega$, as the problem presented is lightly damped, therefore the real part of the eigenvalues is much smaller than the imaginary one.

We seek a parametrisation of the invariant manifold with respect to new normal coordinates $\bfz$ in the form:
\begin{equation}\label{eq:polyrep_map}
    \bfW(\bfz) 
    = \sum_{p=1}^{o}\sum_{k=1}^{m_p} \bfW^{(p,k)} \bfz^{\alphavec (p,k)},
\end{equation}
on the manifold, the dynamics is also described by a polynomial expansion as:
\begin{equation}\label{eq:polyrep_dyn}
    \dot{\bfz} = \bff(\bfz) 
    = \sum_{p=1}^{o}\sum_{k=1}^{m_p} \bff^{(p,k)} \bfz^{\alphavec (p,k)},
\end{equation}
 
The monomial associated to $\alphavec (p,k)$ simply reads $\bfz^{\alphavec (p,k)} = z_1^{\alpha_1}z_2^{\alpha_2} \hdots z_{d}^{\alpha_{d}}$. 

In the context of a non-autonomous problem, the normal coordinates are of two different types: the coordinates parametrising the autonomous manifold $\bar{\bfz}$, and those accounting for the time dependence of the invariant manifold $\tilde{\bfz}$:
\begin{equation}\label{eq:zinormalcoor}
    \bfz = 
    \begin{bmatrix}
        \bar{\bfz}
        \\
        \tilde{\bfz},
    \end{bmatrix}.
\end{equation}

In particular, the first set of variables $\bar{\bfz}$ represent an identity-tangent transformation of the modal coordinates of the selected master subspace. Typically, the forcing is resonant with only one mode, so a single complex conjugate pair of $\bar{\bfz}$ coordinates is enough to reproduce the system's behaviour. Regarding the second set of variables $ \tilde{\bfz}$, they are book-keeping coordinates used in the implementation to automatically expand in the $\varepsilon$ parameter, thus avoiding first-order truncation on the amplitude of the forcing~\cite{JAIN2021How,opreni22high}, see~\cite{vizza2023superharm} for a more detailed presentation. They are defined as:
\begin{align}
\label{eq:tildez}
    \tilde{\bfz} = \begin{bmatrix}
        \varepsilon \e^{+\iu \omega t}
        \\
        \varepsilon \e^{-\iu \omega t}
    \end{bmatrix}.
\end{align}
When the change of coordinates $\bfW(\bfz)$ is inserted into the original equation in place of $\bfy$, the invariance equation is obtained:
\begin{align}
    \label{eq:invariance}
    {\bfB}
        \nabla \bfW(\bfz) \bff(\bfz)
    =
    {\bfA}
    \bfW(\bfz)
    +
    \bsnQ(\bfW(\bfz), \bfW(\bfz))
    +
    \bsnT(\bfW(\bfz), \bfW(\bfz), \bfW(\bfz))
    +
    \begin{bmatrix}
        \bfC_+ & \bfC_-
    \end{bmatrix} 
    \tilde{\bfz},
\end{align}
where the first term comes from the chain rule applied to $\dot{\bfy}$.



The invariance equation is solved recursively, order by order. At order one, the linear problem is retrieved, while higher orders from $p=2$ give rise to the so-called homological equation of order $p$. These equations can be written at the level of an arbitrary monomial, resulting in an underdetermined, ill-conditioned problem. The under-determinacy comes from the fact that two unknowns are input, the nonlinear mapping and the reduced dynamics. This problem is cured by introducing different styles of solution, the most important ones being the graph style and the normal form style~\cite{Haro}. The ill-conditioning comes from the nonlinear resonance relationships. In our implementation of the DPIM, this is cured by augmenting the size of the problems to solve with a bordering technique\cite{vizza21high,opreni22high,vizza2023superharm}.

At the end of the process, the reduction method provides the nonlinear mapping~\eqref{eq:polyrep_map} and the reduced dynamics~\eqref{eq:polyrep_dyn} at arbitrary order. To study the convergence of the ROM, an expansion order needs to be selected. Since the two different parts of the master coordinate $\bfz$ in Eq.~\eqref{eq:zinormalcoor} have different meanings, the truncation of the DPIM ROM will be hereafter referred to as $\mathcal{O}(p,q)$. The first integer $p$ refers to the maximal order of the $\bfz$ coordinate, while $q$ refers to the maximal order of $\tilde{\bfz}$. As a matter of fact, $q$ directly represents the order $\varepsilon^q$ of the truncation of the non-autonomous part. In general, the order $q$ needed to reach convergence is much smaller than that needed for $p$, see the numerical results in Section~\ref{sec:results} and those reported in previous works ~\cite{vizza2023superharm}

Lastly, note that at the end of the procedure, a minimal ROM described by Eq.~\eqref{eq:polyrep_dyn} is obtained. The size of this ROM coincides with the size of $\bfz$, which, for a forcing resonant with a single mode of the system, will be equal to 4. However, since the coordinates $\tilde{\bfz}$ are known from Eq.~\eqref{eq:tildez}, they do not constitute free variables hence the effective size of the ROM is equal to 2. Moreover, different values of the forcing can be simulated from a single run of the DPIM procedure by simply varying the $\varepsilon$ parameter.




\section{Simplified approach for large parallel plates}
\label{sec:PPF}

A model problem  which provides an excellent benchmark for the proposed formulation is represented by a long clamped-clamped beam facing a parallel rigid electrode of the 
same length with initial gap $g$, small with respect to the length.
In this case a very accurate analytical approximation of the electrostatic forces $\bfb$
is given by the {\it parallel plates} formula which assumes an electric vector field that is orthogonal to the mid-plane of the plates, and thus colinear with the $\bfe_2$ vector in our notations, yielding:
\begin{align}
\label{eq:anes}
\bfb(\bfx)=
\frac{1}{2}\veps_0\frac{(V_{\tDC}+V_{\tAC}\sin\omega t)^2}{(g-u_2(\bfx))^2}
\bfe_2,
\end{align}
where $\bfx$ stands for the point position 
in the initial configuration and the force is exerted only on the upper surface $\tU_0$ facing the electrode.
This approximation ignores fringing fields and horizontal components of the force, but is known to be extremely accurate for long beams and small gaps.
In this case, a simplified formulation of the {\sc mgf} can be obtained by using the analytical expression for $\bfb$, yielding:
%
\begin{equation}
\label{eq:AF}
\int_{\Omega_0}\rho_0\ddot{\bfu}\cdot\tilbfu\,\dOmega_0 + 
\int_{\Omega_0} \bfs[\bfu]:(\bfeps[\tilbfu] + \bfens[\bfu,\tilbfu])\,\dOmega_0 
=
\frac{\veps_0}{2}\int_{\tU_0}\frac{(V_{\tDC}+V_{\tAC}\sin\omega t)^2}{(g-u_2(\bfx))^2}\tilde{u}_2\,\dS_0,
\quad \forall \tilbfu \in \cC_s(0).
\end{equation}
This model will be referred to as the parallel-plate formulation ({\sc ppf}) in the rest of the paper. The analytical expression of $\bfb$ substantially simplifies the problem such that direct time integration routines can be directly applied to Eq.~\eqref{eq:AF}. The {\sc ppf} will thus be used in the next developments to provide time-domain full-order reference solutions.
However, Eq.~\eqref{eq:AF} still
contains non-polynomial nonlinearities, which represents a major obstacle for several solution schemes and for the DPIM in particular. In analytical approaches, like the multiple scales technique used for instance in \cite{younis03,younisbook}, two solutions are adopted.
The most direct approach consists in expanding $1/(g-u_2(\bfx))^2$ 
in Taylor series. However, the series converges very slowly and a very large number of terms
is required when $u_2$ reaches significant fractions of the gap.
As an alternative {\it exact} approach which is also strictly connected with the general formulation developed in the previous section, 
one can introduce  an auxiliary field:
\begin{align}
\label{eq:psi}
\psi(\bfx)=\frac{V_{\tDC}+V_{\tAC}\sin\omega t}{g-u_2(\bfx)}.
\end{align}
Remark that $\psi$ has the meaning of a gradient of the potential, so that a notation similar to $\bfpsi$ in the general formulation of the previous sections is adopted.
The power of electric forces are quadratic in $\psi$ and
Eq.~\eqref{eq:AF} becomes:
\begin{align}
\label{eq:MAFmech}
\int_{\Omega_0}\rho_0\ddot{\bfu}\cdot\tilbfu\,\dOmega_0 + 
\int_{\Omega_0} \bfs[\bfu]:(\bfeps[\tilbfu] + \bfens[\bfu,\tilbfu])\,\dOmega_0 
=\frac{\veps_0}{2}\int_{\tU_0} \psi^2(\bfx)\tilu_2\dS_0, 
\quad \forall \tilbfu \in \cC_s(0).
\end{align}
Moreover, the strong form \eqref{eq:psi} is replaced with a weak form
of the {\it compatibility equation}:
\begin{align}
\label{eq:MAFcomp}
\int_{\tU_0} \psi(\bfx)(g-u_2(\bfx))\tilpsi(\bfx)\dS_0
=(V_{\tDC}+V_{\tAC}\sin\omega t)\int_{\tU_0}\tilpsi(\bfx)\dS_0,
\quad \forall\,\tilpsi(\bfx)\in \cC_\psi,
\end{align}
where $\tilpsi$ and $\psi$ are chosen in the same space $\cC_\psi$.
In Section \ref{sec:results}, the results presented 
are generated by an implementation where $\cC_\psi$ is the space of piecewise constant
shape functions.

Eqs.~\eqref{eq:MAFmech}-\eqref{eq:MAFcomp} represent the mixed semi-analytical formulation for parallel plates problems only, and is referred to as 
{\sc mppf}. It can be discretised following the same steps as for the {\sc mgf}
and then reduced using the DPIM.
Details are collected in Appendix~\ref{app:mafdiscr}.
The {\sc mppf} model will be used to validate numerical results in
Section~\ref{sec:results}.

\section{Numerical results}
\label{sec:results}


In this section, we collect some numerical tests intending to validate the proposed formulations.

First, a case of static pull-in instability is addressed in Section~\ref{sec:pullin}
using the mixed electrostatic formulation introduced in Section~\ref{sec:elecinitial}
and a graphical representation of the new associated fields is also provided for the sake of illustration. 
Next, in Section~\ref{sec:compaformulationROM}
the different ROM strategies are applied to the case of a Parallel Plate setting, where  
a clamped-clamped beam faces a parallel rigid electrode of the same length.
Despite its simple geometrical layout, this configuration activates several nonlinear
phenomena when computing either  primary or secondary resonances, and hence represents
a challenging testbed.

\subsection{Pull-in verification}
\label{sec:pullin}

To validate the proposed mixed general formulation {\sc mgf} for the electrostatic field,
we examine the static pull-in phenomenon that occurs when an electrostatically actuated system reaches a critical voltage, causing its movable component to abruptly snap into contact with a fixed electrode due to electrostatic forces.

\begin{figure}[ht]
\centering
\includegraphics[width = .7\linewidth]{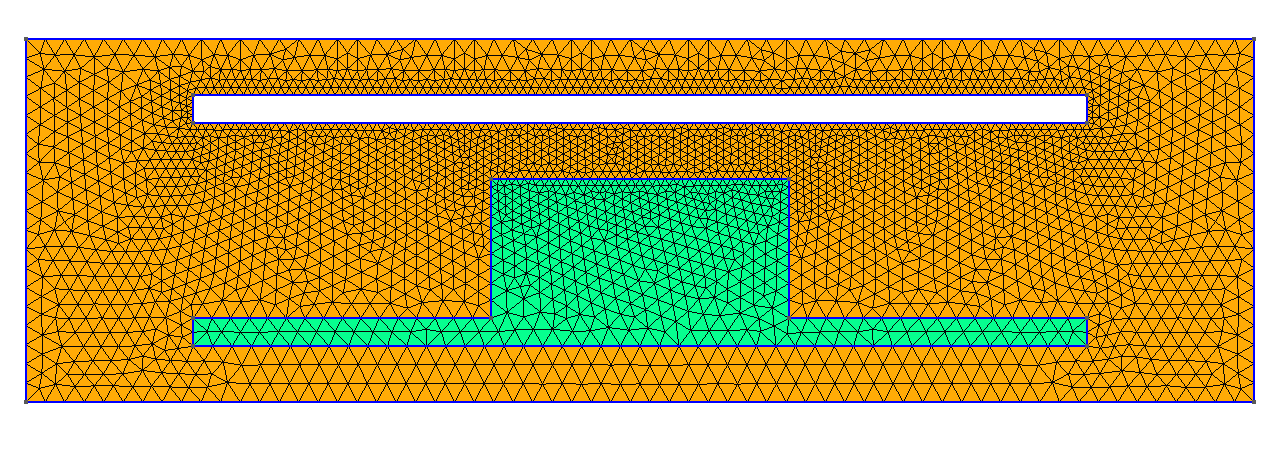}
\caption{Static pull-in verification. Green mesh: deformable structure $\Omega_0$;
orange mesh: air domain $\tA_0$.
The lower deformable mass is anchored by two thin beams on either side to model linear spring constraints.} 
\label{fig:pullin_mesh}
\end{figure}

A classical result concerns parallel plates in which a rigid plate, with given voltage $V_\tDC$ and constrained by linear springs, faces a parallel fixed electrode with zero voltage. 
When $V_\tDC$ increases quasi-statically, the former plate stably approaches the electrode until the displacement is equal to one-third of the gap.
Above this value, the system reaches the pull-in instability, and a further increase of the voltage induces a dynamic evolution where the movable plate will snap to the electrode.
To reproduce such a set-up numerically, the case shown in Fig.~\ref{fig:pullin_mesh} is considered.
The fixed upper electrode is shown in white, while the movable part, in green, consists of a central stiff mass, constrained by two thin beams on either side.
The nonlinear mechanical effects in these beams are deliberately disabled to simulate linear springs. The problem is two-dimensional, and the mass has a length of $6$\micr and a height of $3$\micr, while the beam thickness is $0.5$\micr for a length of $5$\micr. The gap $g$ between the electrode and the mass is $1$\micr.
The elastic constants are set to $E = 154$\,GPa and $\rho=2330$\,kg/m$^3$.

\begin{figure}[ht]
\centering
\includegraphics[width = .7\linewidth]{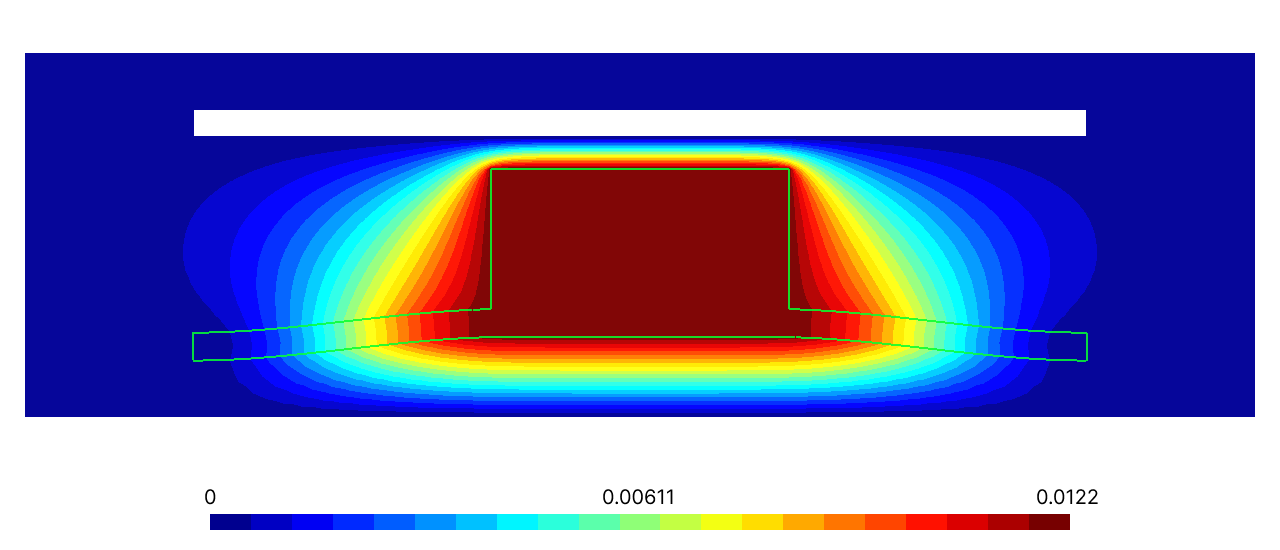}
\caption{Static pull-in verification. Static analysis with $V_\tDC=10$\,V.
Vertical displacement $u_{02}$ in\micr in the solid and its extension $s_{02}$ to the air domain. 
The two fields are continuous on the solid surface. The field is plotted in the deformed configuration for graphical purposes, but the formulation is fully Lagrangian.
}\label{fig:pullin_u0}
\end{figure}

\begin{figure}[p]
\captionsetup[sub]{font=normalsize,labelfont={bf,sf}}
\begin{subfigure}{1.0\columnwidth}
\centering
\includegraphics[width = 0.8\linewidth]{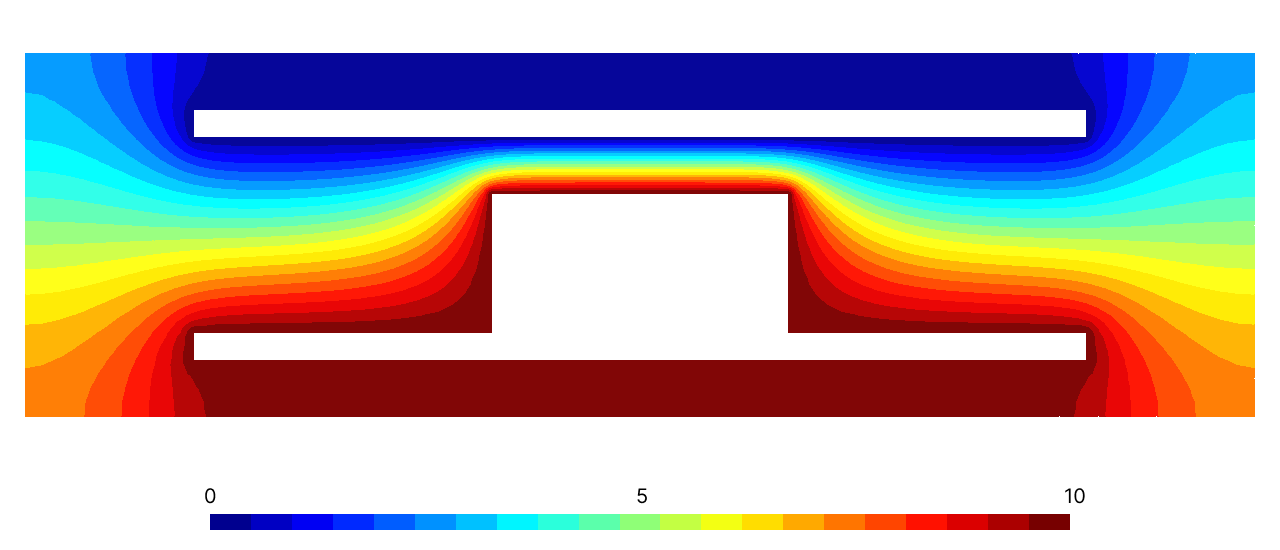}
\caption{$\varphi_{0}$}
\label{fig:pullin_phi}
\end{subfigure}
\begin{subfigure}{1.0\columnwidth}
\centering
\includegraphics[width = 0.8\linewidth]{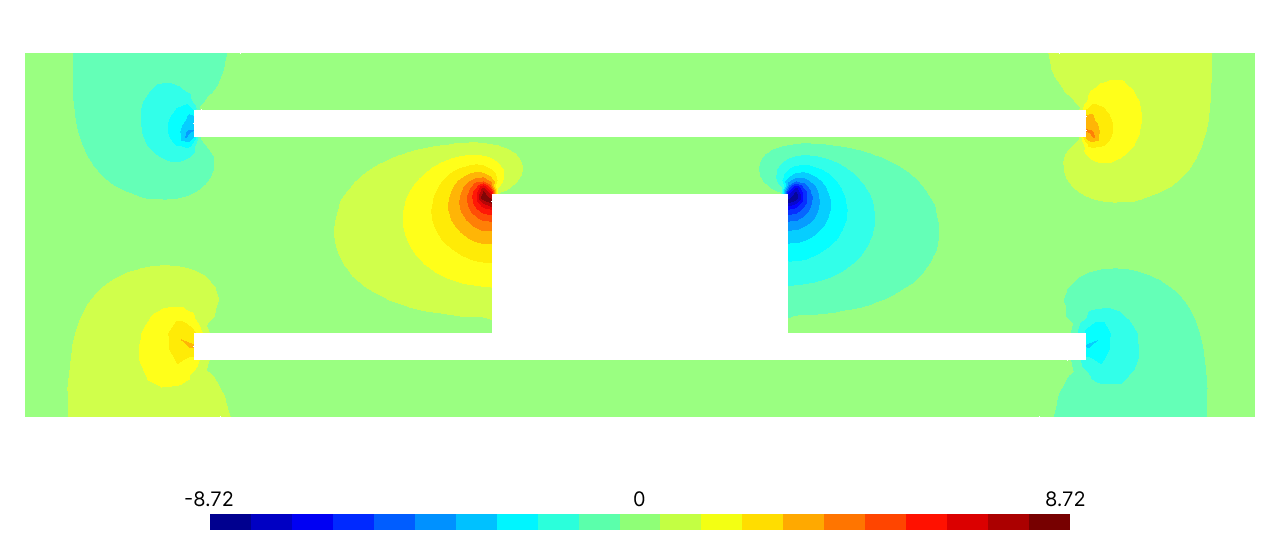}
\caption{$\psi_{01}$}
\label{fig:pullin_psi1}
\end{subfigure}
\begin{subfigure}{1.0\columnwidth}
\centering
\includegraphics[width = 0.8\linewidth]{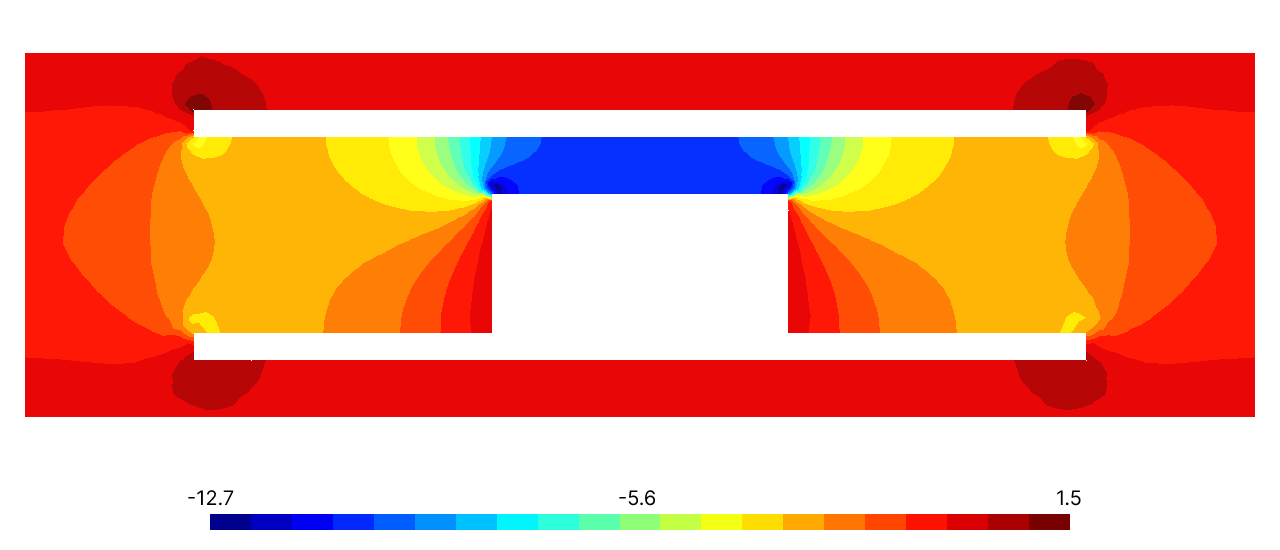}
\caption{$\psi_{02}$}
\label{fig:pullin_psi2}
\end{subfigure}
\caption{Static pull-in verification. Static analysis with $V_\tDC=10$\,V.
(a) potential $\varphi_0$. The values of $V_\tDC$ assigned on the electrode and the body, respectively $0\,$V and $10\,$V, correspond to $\varphi_0$ on the surface of the bodies. Remark the linear distribution between the plates;
(b) Horizontal component $\psi_{01}$ of the Eulerian gradient showing singularities near the edges of the electrodes;
(c) Horizontal component $\psi_{02}$ of the Eulerian gradient which is nearly constant within the parallel plates,  corresponding to the presence of a uniform electric field.}
\label{fig:pullin_qoi}
\end{figure}

Figures~\ref{fig:pullin_u0} and \ref{fig:pullin_qoi} present the solution fields
for the particular case of $V_\tDC=10$\,V.
In particular, Figure~\ref{fig:pullin_u0}
is the contour plot of the vertical displacement of the structure $u_{02}$, and its extension $s_{02}$ to the air domain $\tA_0$.
The overall field is continuous and vanishes on $S_{\infty}$ and on the surface of the electrode.

The plots of Figures~\ref{fig:pullin_phi}-\ref{fig:pullin_psi2} correspond to the potential $\vphi_0$ and the two components of the Eulerian gradient $\bfpsi_{0}$, 
Eq.~\eqref{eq:new}, always for $V_\tDC=10$\,V. 
It can be appreciated that the electric field is almost constant and vertical between the plates, which justifies the comparison with parallel-plate formulas for the pull-in value.

\begin{figure}[ht]
\centering
\includegraphics[width = .7\linewidth]{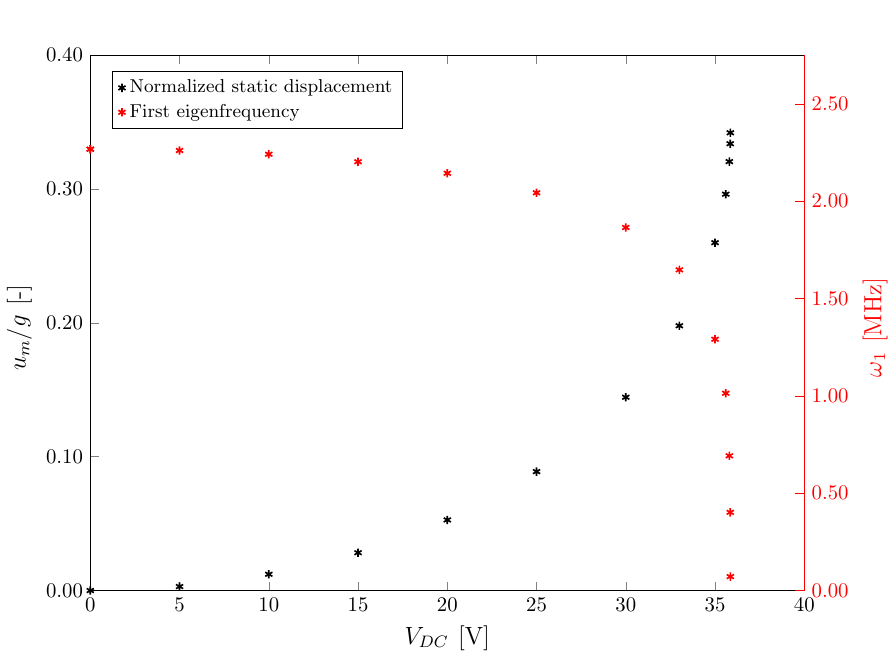}
\caption{Static pull-in verification. Evolution of the relative displacement 
$u_m/g$ (left axis) of the lower rigid mass and of the first eigenfrequency (right axis) 
at different actuating voltages. As a consequence of the pull-in phenomenon, when the displacement approaches one-third of the initial gap, even a small voltage increase causes a large change in the displacement and the first eigenvalue vanishes.}
\label{fig:pullin_u_omega}
\end{figure}

Next, the voltage bias between the electrode and the top side of the mass is raised, and the static deformation at different voltage levels is computed. Figure~\ref{fig:pullin_u_omega} reports the evolution of the maximum displacement and the first eigenvalue $\omega_1$ as the potential varies.
The plot clearly shows that the mass displacement increases with the actuating voltage in a stable manner; nevertheless, as the midspan displacement $u_m$ approaches the critical value of one-third of the initial gap, even small increments in the actuating voltage result in significant changes in deformation, following the typical behaviour observed in pull-in phenomena.
Similarly, the first eigenfrequency gradually decreases due to the softening influence of the potential difference, then rapidly converges to $0$ as the critical voltage threshold is approached, confirming that the instability phenomenon is correctly captured.

\begin{figure}[p]
\captionsetup[sub]{font=normalsize,labelfont={bf,sf}}
\begin{subfigure}{1.0\columnwidth}
    \centering
    \includegraphics[width = 0.7\linewidth]{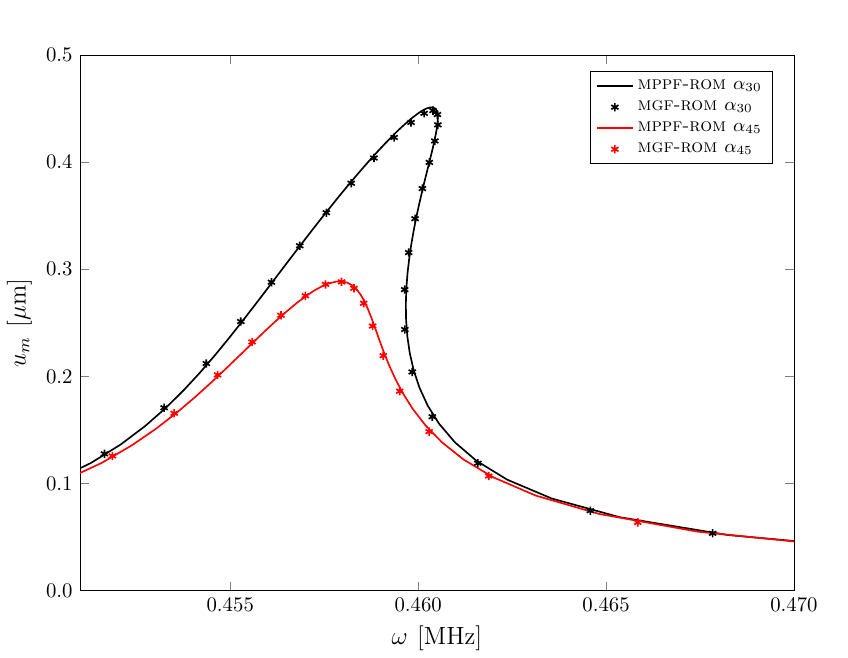}
    \caption{ }
    \label{fig:T3ANhard}
\end{subfigure}
\begin{subfigure}{1.0\columnwidth}
    \centering
    \includegraphics[width = 0.7\linewidth]{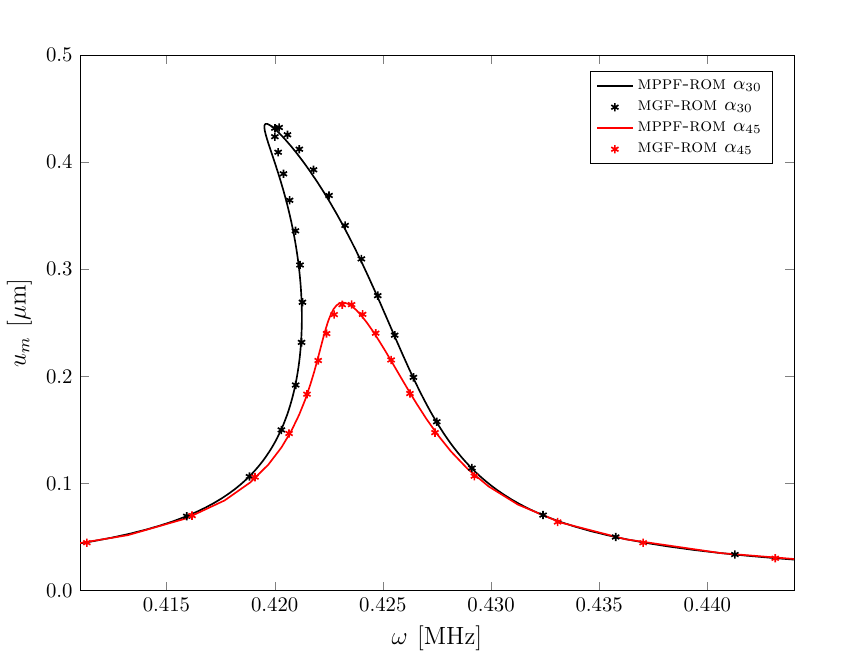}
    \caption{ }
    \label{fig:T3ANsoft}
\end{subfigure}
\caption{Parallel plates. Comparison between the responses of the Reduced Order Models {\sc mppf-rom} and {\sc mgf-rom}. Frequency response curves for the mid-span displacement $u_m$. (a) Hardening case: $V_\tDC=1.2$\,V, $V_\tAC=0.2$\,V,
with two different levels of dissipation $\alpha_{30}=3\,10^{-3}s^{-1}$ and  $\alpha_{45}=4.5\,10^{-3}s^{-1}$. (b) Softening case $V_\tDC=4.0$\,V, $V_\tAC=0.045$\,V
with the same dissipation levels. Increasing the $V_\tDC$ value, the negative spring effect due to the electrostatic field dominates and the overall beam response turns from hardening to softening.}
\label{fig:T3ANhard+soft}
\end{figure}

\subsection{Reduced order modelling for parallel pates}
\label{sec:numpp}

We now turn to the analysis of the dynamical response of a clamped-clamped microbeam of length $L =510$\micr, 
thickness $h=1.5$\micr and air-gap $g =1.18$\micr
with respect to a parallel electrode that has the same length.
The dimensions correspond to those of realistic microstructures,
see for instance \cite{younis03,younis05,younisbook}, and 
are chosen to trigger the wealth of nonlinear effects pointed out in what follows.
The microbeam is excited by a voltage bias of the form $V_\tDC + V_\tAC \sin(\omega t)$.
We also assume that the beam is subjected to an axial force of 0.0009 N.
Elastic constants are $E =154$\,GPa and $\rho=2330$\,kg/m$^3$
while the Poisson coefficient is set to $0$ to simulate a large thin plate.

In the following section, the various formulations proposed in this paper are used to generate the Frequency Response Curves (FRCs) 
for the mid-span displacement $u_m$ of the beam, i.e. the amplitude of the oscillation
of the center point when the frequency $\omega$ of the $V_\tAC$
input signal spans a specific range of frequencies.

In all the applications, a single master mode is selected, excluding a priori possible phenomena of internal resonance. The master mode corresponds to the first bending mode of the clamped-clamped beam interacting with the surrounding electric field.
This implies that the size of the normal coordinates vector
$\bfz$ in the reduced order model $\dot{\bfz} = \bff(\bfz)$ in 
Eq.\eqref{eq:polyrep_dyn}
is 4. The first two coordinates, $z_1$ and $z_2=\bar{z}_1$, correspond to the normal coordinates representing the master mode, while the last two, $(z_3,z_4)$, are the non-autonomous variables representing the harmonic forcing, such that $z_3=\e^{\iu \omega t}$ and $z_4=\e^{-\iu \omega t}$. 

Both primary and secondary resonances will be investigated in the following sections.
For convenience, we gather in Table~\ref{tab:cases} the values utilized 
in the different analyses for
$V_{\tDC}$, $V_{\tAC}$ and for the actuating frequency $\omega$.

\begin{table}[h]
\centering
\begin{tabularx}{0.8\textwidth}{p{5cm}>{\centering\arraybackslash}X>{\centering\arraybackslash}X>{\centering\arraybackslash}X}
\hline
 & $V_{\tDC}$ [V] & $V_{\tAC}$ [V] & $\omega$ [Hz] \\
\hline
Primary resonance - hardening & 1.2 & 0.2 & $\simeq\omega_1$ \\
Primary resonance - softening & 4.0 & 0.045 & $\simeq\omega_1$ \\
Superharmonic resonance & 1.7 & 0.7 & $\simeq 0.5\,\omega_1$ \\
Subharmonic resonance & 1.7 & 0.7 & $\simeq 2\,\omega_1$ \\
\hline
\end{tabularx}
\caption{Summary of the test cases performed}
\label{tab:cases}
\end{table}

\subsubsection{Comparison between the reduced order models of the general 
{\sc mgf} and parallel-plate {\sc mppf} formulations}
\label{sec:compaformulationROM}

This section aims to validate the parallel plate assumptions by comparing the results provided by the {\sc mgf} with those of the {\sc mppf}, with the ultimate goal of using
in the next section the more agile and cheaper {\sc mppf} to provide full-order reference solutions. 
As a consequence, we start by comparing the outcomes of the ROMs obtained from the two formulations. Since the DPIM technique has already been shown to deliver very accurate results for purely mechanical vibrating systems, see {\it e.g.} \cite{vizza21high,vizza2023superharm,opreni22high}, we will use this strategy to validate the fact that {\sc mppf} provides trustworthy solutions, in line with numerous previous studies showing the accuracy of the parallel plate analytical formula.

\begin{figure}[ht!]
\captionsetup[sub]{font=normalsize,labelfont={bf,sf}}
\begin{subfigure}{1.\columnwidth}
    \centering
    \includegraphics[width = .7\linewidth]{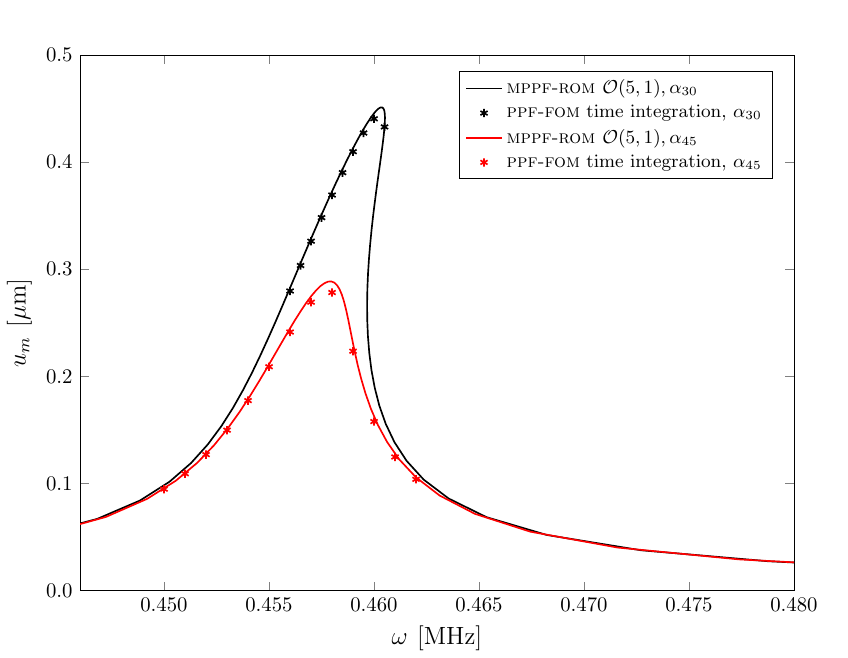}
    \caption{ }
    \label{fig:timehard}
\end{subfigure}
\begin{subfigure}{1.\columnwidth}
    \centering
    \includegraphics[width = .7\linewidth]{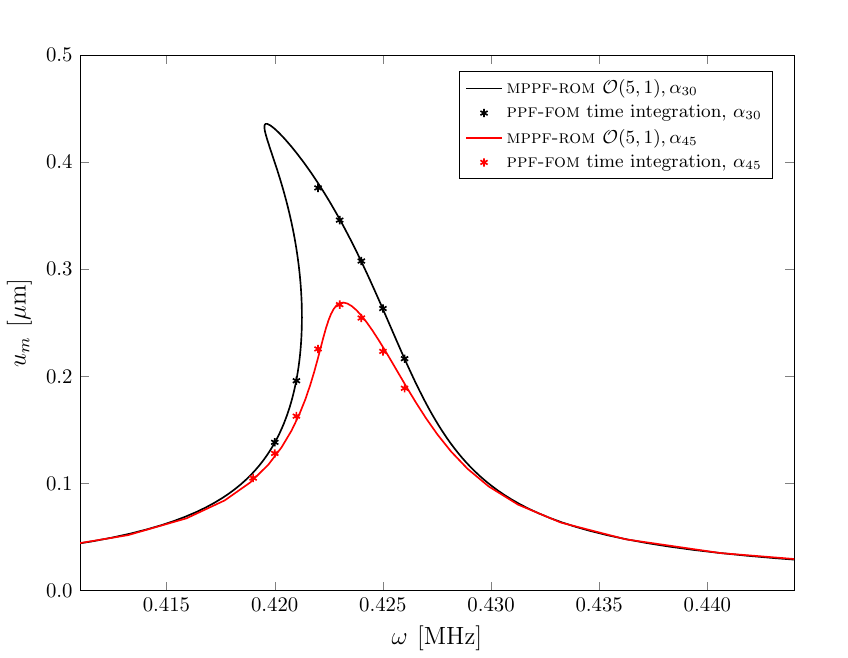}
    \caption{ }
    \label{fig:timesoft}    
\end{subfigure}
\caption{
Parallel plates. 
Comparison of the full Full Order Model {\sc ppf}, time-integrated,
with the Reduced Model {\sc mppf-rom}, solved with a  numerical continuation procedure. 
The same actuation levels as in Section \ref{sec:compaformulationROM} are
utilized.
(a) $V_\tDC=1.2$\,V, $V_\tAC=0.2$\,V. (b) $V_\tDC=4$\,V, $V_\tAC=0.045$\,V.
The two dissipation levels are:
$\alpha_{45}=4.5\,10^{-3}s^{-1}$, $\alpha_{30}=3.0\,10^{-3}s^{-1}$}
\end{figure}

Two cases of actuation are selected to draw out this comparison. First, a small value of $V_\tDC$ is imposed to get a hardening behaviour for the vibration of the fundamental mode. Next, $V_\tDC$ is increased and the negative spring effect due to the electrostatic field becomes dominant, turning the overall beam response from hardening to softening, a well-known effect in electrostatic MEMS.
For each formulation, {\sc mgf} and {\sc mppf}, a ROM has been derived with the DPIM up to order $\mathcal{O}(5,1)$ using the complex normal form style. 
The ROM is then solved by means of a continuation technique implemented in the software  Matcont~\cite{dhooge2004matcont} to obtain the frequency-response curve (FRC) in the vicinity of the fundamental eigenfrequency. 


Figure \ref{fig:T3ANhard+soft} collects the results for two different sets of parameters.
The case of Figure \ref{fig:T3ANhard} corresponds to $V_\tDC=1.2$\,V, $V_\tAC=0.2$\,V 
and $\alpha_{30}=3\,10^{-3}s^{-1}$, $\alpha_{45}=4.5\,10^{-3}s^{-1}$.
Here the $V_\tDC$ bias is small and the beam retains the classical hardening behavior 
of a purely mechanical clamped-clamped beam.
In Figure \ref{fig:T3ANsoft}, on the contrary, we set
$V_\tDC=4$\,V and $V_\tAC=0.045$\,V, a situation much closer to pull-in.

Note that, for the sake of simplicity, the stability analysis is not reported in the figures. Since the cases under study refer to classical cases in nonlinear vibrations, the stability of the different branches of solution can be easily inferred from the results, and are correctly predicted by the ROMs. The level of actuation has been set such that the maximum vibration amplitude reported in Fig.~\ref{fig:T3ANhard+soft} corresponds to 1/3 of the thickness and more than 1/3 of the gap $g$, hence giving rise to important nonlinearities with respect to the setup.

In the numerical results, two effects should be noticed.  First, there is a shift of the overall FRC towards 
smaller frequencies. Indeed,
increasing the $V_\tDC$ value the negative spring effect due to the electrostatic field overcomes the expected mechanical hardening.

Secondly, the nonlinear response of the beam is now dominated by electrostatic forces and becomes softening.
In general, the agreement between the two different approaches is very satisfactory, undoubtedly showing that for this application and the range of parameters used, the simplification of the electric field is perfectly well justified. Consequently, the {\sc mppf} formulation will be used in the rest of the paper to produce the FOM solutions.

\subsubsection{Comparison between the {\sc mppf} full-order model and the ROM}
\label{sec:FOMROM}

Since the reduced models obtained with the {\sc mgf} and {\sc mppf} are fully equivalent, as shown in the previous section,
we now compare the results obtained from the {\sc mppf-rom}
with those stemming from a time integration approach, directly applied to the full-order {\sc mppf} model, Eq.\eqref{eq:AF}, which will be used as FOM reference solution.

In particular, we analyze the same cases of hardening and softening behaviour as in the previous section.
In Figure~\ref{fig:timehard}, the voltages are set to $V_{\tDC} = 1.2$\,V and $V_{\tAC} = 0.2$\,V,
while in Figure\ref{fig:timesoft}, $V_{\tDC} = 4$\,V and $V_{\tAC} = 0.045$\,V. The structural losses are modelled with a mass-proportional Rayleigh damping $\bfC = \alpha \bfM$, and the Frequency Response Curves are computed for two different dissipation values, respectively $\alpha_{45}=4.5\,10^{-3}s^{-1}$ and $\alpha_{30}=3\,10^{-3}s^{-1}$.

The FRCs of the {\sc mppf-rom} 
is obtained with the continuation package Matcont,
while the FOM is integrated in time with a nonlinear implicit Newmark scheme, with $\beta=1/4,\gamma=1/2$.
In the latter case, each cycle is partitioned in $120$ steps and $3Q$ cycles are simulated for each value of $\omega$ always using a homogeneous state as an initial condition for the analysis, where $Q$ is the quality factor. This time integration allows to follow only the stable branches, while the continuation technique can reproduce the whole Frequency Response Curve.

In all four cases, the ROM and FOM results are perfectly aligned, providing additional confirmation of the extreme accuracy of the reduced approach.


\subsubsection{Convergence study}
\label{sec:convergence}

\begin{figure}[ht!]
\captionsetup[sub]{font=normalsize,labelfont={bf,sf}}
\begin{subfigure}{1.\columnwidth}
\centering
\includegraphics[width = .7\linewidth]{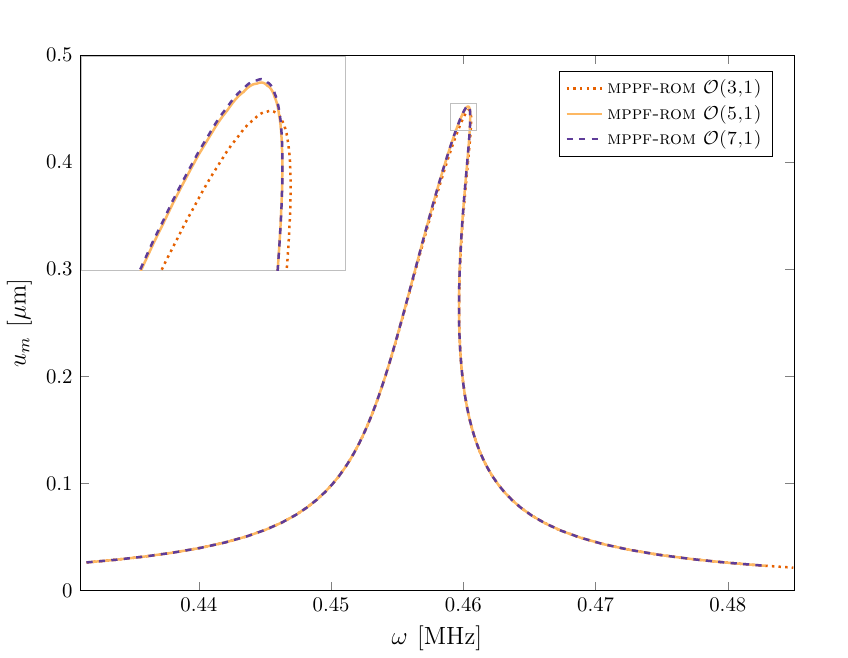}
\caption{}
\label{fig:hard_alpha3}
\end{subfigure}
\begin{subfigure}{1.\columnwidth}
\centering
\includegraphics[width = .7\linewidth]{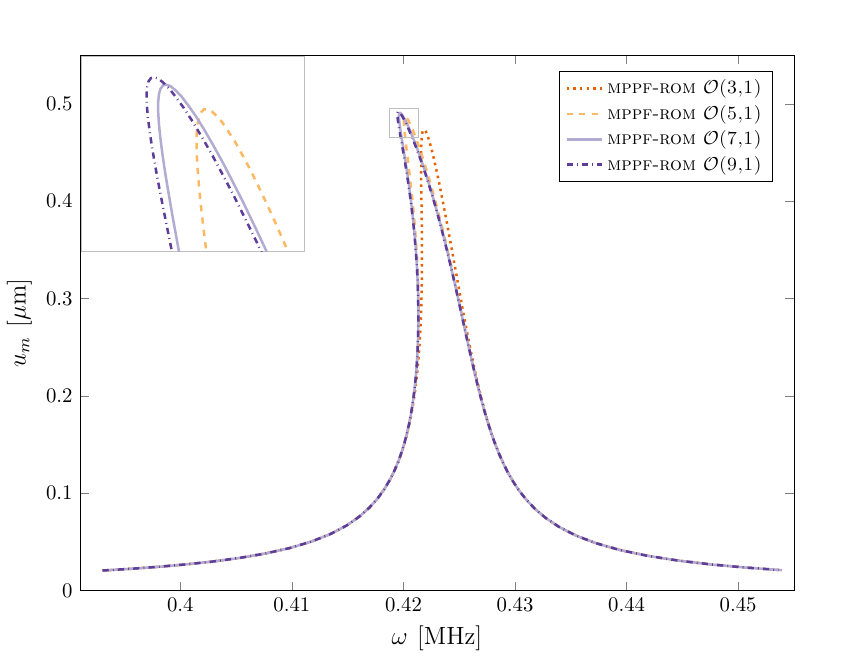}
\caption{}
\label{fig:soft_alpha3}
\end{subfigure}
\caption{
Parallel plates. 
Convergence analysis of the ROM for the primary resonance. (a): hardening behaviour with low damping. $V_{\tDC} = 1.2$\,V, $V_{\tAC} = 0.2$\,V, $\alpha=3\,10^{-3}$\,s$^{-1}$;
(b) softening behaviour with low damping. $V_{\tDC} = 4.0$\,V, $V_{\tAC} = 0.045$\,V, $\alpha=3\,10^{-3}$\,s$^{-1}$. }
\end{figure}

Having showcased the {\sc mppf-rom} capability to yield accurate results, aligning closely with the solution provided by the full model, this section aims to assess the convergence rate of the DPIM towards the actual solution while increasing the expansion order.
As the analyses described in the previous section already made apparent, the DPIM exhibits notably rapid convergence in the case of primary resonance.
Specifically, in the hardening test case, convergence is achieved as early as the $\mathcal{O}(5,1)$ expansion order, as shown in Figure \ref{fig:hard_alpha3}.
However, when softening effects are more prominent and the pull-in instability is approached, the rate of convergence tends to slow down: Figure \ref{fig:soft_alpha3} underscores that there remains a noticeable difference between order $\mathcal{O}(5,1)$ and $\mathcal{O}(7,1)$ expansions, indicating the necessity of high order expansions 
near pull-in. Pushing the DPIM order further only yields minimal improvement, as the order $\mathcal{O}(9,1)$ expansion almost overlaps with the $\mathcal{O}(7,1)$ one. Concerning the convergence of the non-autonomous development, which is set to 1 in the previous study, it has been numerically verified (and not shown here for the sake of brevity) that increasing the order does not change the FRC which is indeed converged with an order 1 on the forcing term. This result is in line with previous studies on primary resonances, which were limited to this order of accuracy on the forcing term and showed a good convergence rate, see for example the discussions and examples reported elsewhere~\cite{vizza2023superharm,opreni22high,li2021periodic}.



\subsubsection{Secondary resonances}
\label{sec:subsuper}

In this final Section, the ROM is used to derive predictions on the nonlinear vibrating behaviour of the coupled system when the forcing frequency $\omega$ is in the vicinity of the secondary resonances, corresponding to $\omega \simeq \omega_1/2$, and $\omega \simeq 2\omega_1$ (with $\omega_1$ the linear frequency of the fundamental bending mode), underlining its excellent predictive capabilities in different dynamical scenario.

\begin{figure}[ht!]
\centering
\includegraphics[width = .7\linewidth]{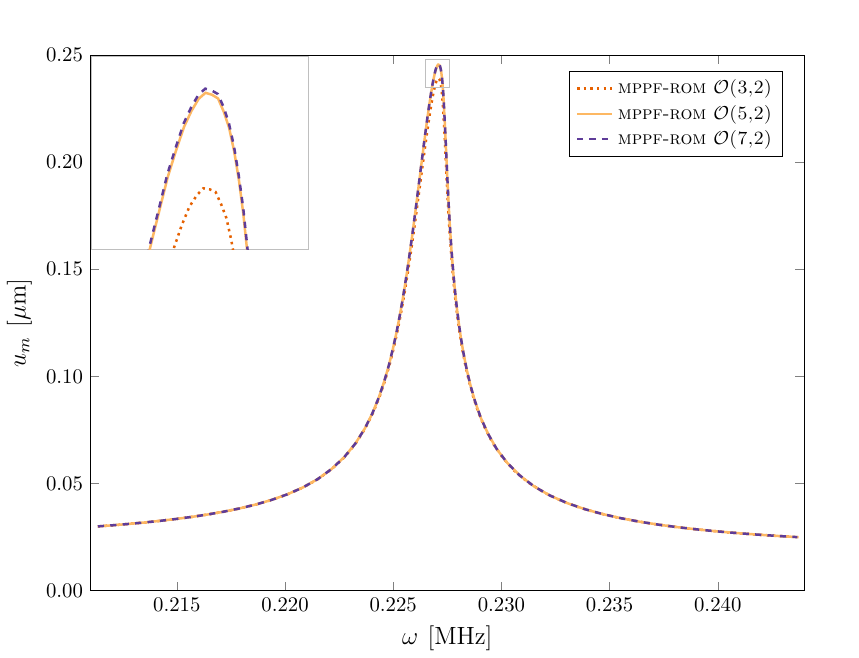}
\caption{Parallel plates. 
Superharmonic resonance with $\omega \simeq \omega_1/2$. $V_{\tDC} = 1.7$\,V, $V_{\tAC} = 0.7$\,V, and mass-proportional damping with $\alpha=3\,10^{-3}s^{-1}$.}
\label{fig:superharm_alpha3}
\end{figure}


The first case under study is that of the superharmonic resonance, when the frequency of the actuation voltage approaches $\omega_1/2$.
As shown in \cite{vizza2023superharm}, a first-order expansion on the non-autonomous term is not able to cope with such superharmonic resonance, and the ROM needs to include a $\varepsilon$ order that is at least equal to the order of the superharmonic considered, 2 in this case. Fig.~\ref{fig:superharm_alpha3} reports the results found when the actuation is in the vicinity of half the eigenfrequency of the fundamental bending mode. Three ROMs with increasing orders are computed to achieve convergence: $\mathcal{O}(3,2)$, $\mathcal{O}(5,2)$ and $\mathcal{O}(7,2)$. Whereas the truncation order  $\mathcal{O}(3,2)$ is not converged, the overlap of the two other curves underlines that an order  $\mathcal{O}(5,2)$ is sufficient to reproduce this nonlinear secondary resonance. While lower orders of the truncation of the non-autonomous term are unable to capture this superharmonic behaviour, it has also been verified numerically that increasing further the order does not change the picture, showing that the results shown in Fig.~\ref{fig:superharm_alpha3} are correctly converged. One can note that such behaviour on secondary resonances has already been reported\cite{Nayfehsub01,younis05}, with a simplified model for the beam and the electrostatic coupling. Here the model computes the structural vibrations with an FE model without kinematical assumptions, and fully takes into account the electromechanical coupling, offering quantitative predictions on the amplitude levels one can reach in such superharmonic regime that could be directly compared {\it e.g.} to experiments.


\begin{figure}[ht!]
\centering
\includegraphics[width = .7\linewidth]{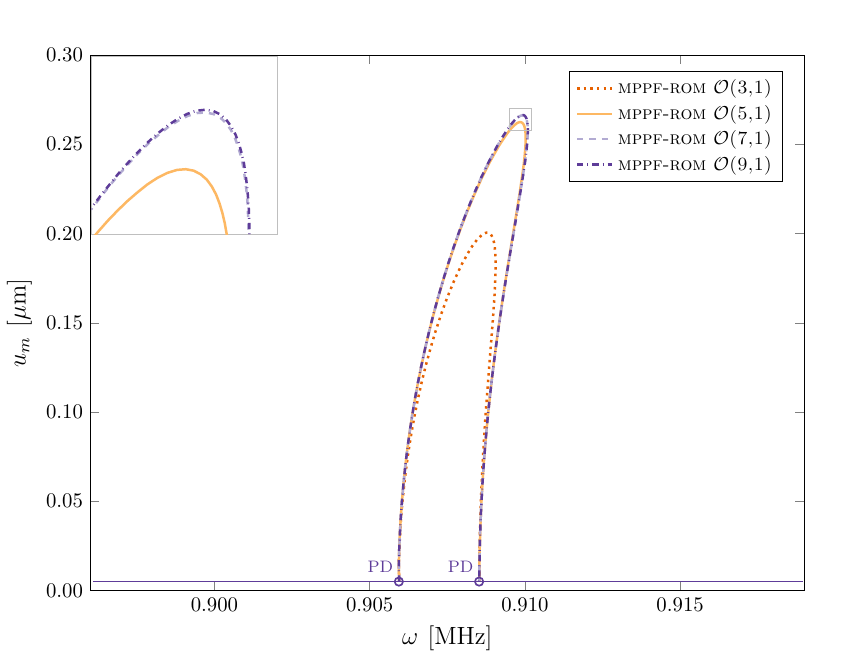}
\caption{Parallel plates. 
Parametric resonance with Rayleigh damping. Convergence of the autonomous part of the expansion. $V_{\tDC} = 1.7$\,V, $V_{\tAC} = 0.7$\,V, $\alpha=-3.81\,10^{-2}s^{-1}$, $\beta=2\,10^{-1}s$}
\label{fig:param1}
\end{figure}

\begin{figure}[ht!]
\centering
\includegraphics[width = .7\linewidth]{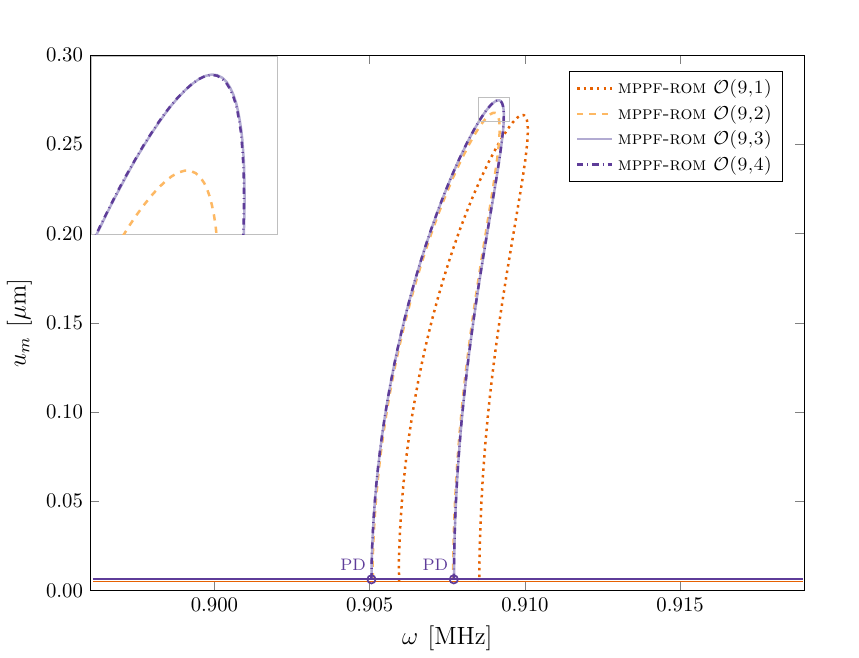}
\caption{Parallel plates. 
Parametric resonance with Rayleigh damping. Convergence of the non-autonomous part of the expansion. $V_{\tDC} = 1.7$\,V, $V_{\tAC} = 0.7$\,V, $\alpha=-3.81\,10^{-2}s^{-1}$, $\beta=2\,10^{-1}s$}
\label{fig:param2}
\end{figure}



\vspace{.5cm}
The last case under study is the secondary resonance occurring when the forcing frequency is in the vicinity of two times the eigenfrequency of the fundamental mode. This case has already been analysed with von K{\'a}rm{\'a}n assumptions for the beam and first-order multiple scales solutions\cite{Nayfehsub01,younis05}, and was termed subharmonic resonance. However, the main driver of the instability is a parametric resonance, as shown in the solvability conditions derived in~\cite{Nayfehsub01}, that contains the same term as the classical parametric resonance~\cite{Nayfeh79,Thomsen}. The resonance curves shown in \cite{Nayfehsub01,younis05} are also typical of a parametric excitation, and not a subharmonic resonance of a single degree of freedom oscillator, which is characterized by a detached solution branch~\cite{Nayfeh79,VOLVERT2021}. As a matter of fact, this resonance scenario can be interpreted in light of the dynamical system equations~\eqref{eq:ODE_original}, and appears to be a subharmonic resonance for the complete system, but a parametric resonance for the mechanical vibratory systems, once the other equations have been condensed into the dynamical problem. For this reason, we will refer to this case hereafter as a parametric excitation.

To give further insight into this parametric resonance, the ROM derived thanks to the parametrisation method can be used as it provides a fruitful physical interpretation of the nonlinear dynamical phenomena at hand. 
Truncating the ROM to order 3 for both the autonomous $(z_1,\bar{z}_1)$ coordinates and the non-autonomous ones $(z_3,z_4)$, the reduction method provides the following reduced dynamics in complex normal form:
 \begin{equation}\label{eq:romparam}
     \dot{z}_1 = \lambda_1 z_1 + f_1 z_1^2 \bar{z}_1  +  f_2 \bar{z}_1z_3 + f_3 z_1 z_3 z_4    + \mathcal{O}(|\bfz|^3),
 \end{equation}
 where $\lambda_1 = -\xi_1\omega_1 + \iu \omega_1\sqrt{1-\xi_1^2}$ is the complex eigenvalue of the fundamental mode, and $\xi_1$ refers to the modal damping ratio. The second monomial corresponds to the trivially resonant cubic term, while the third one represents the parametric excitation term, since $z_3=\e^{\iu\omega t}$, such that $f_2$ is directly proportional to the amplitude of the forcing. Finally, the last monomial $z_1 z_3 z_4$ is a trivially resonant term of higher order in the non-autonomous term and needs an order 2 on the non-autonomous coordinates to be considered.

In Eq.~\eqref{eq:romparam}, the driver of the parametric instability is the second-order term $f_2 \bar{z}_1z_3$ which reads $f_2 \bar{z}_1\e^{\iu\omega t}$.  This monomial is indeed resonant when  $\omega$ is in the vicinity of $2 \omega_1$.  In particular, through a period-doubling bifurcation from the main branch of zero solution, a periodic response occurs in the form:
 \begin{equation}\label{eq:polarrep}
     z_1(t) = \rho e^{\iu \phi} e^{\iu \omega t/2},
 \end{equation}
 where a polar representation of the normal coordinate is introduced with an amplitude $\rho$. In order to show how the ROM predicts the bifurcation points, let us truncate Eq.~\eqref{eq:romparam} to its lowest-order significant terms by neglecting the two monomials of order three with coefficients $f_1$ and $f_3$. Inserting the ansatz~\eqref{eq:polarrep} into the normal form with only the $f_2$ term of the parametric instability, and assuming a steady-state response since we are interested in the solution branches of the permanent regime, one obtains:
 \begin{equation}
     \rho \left( \iu \omega/2 - \lambda_1 - f_2 e^{-2\iu \phi} \right) = 0  + \mathcal{O}(\rho^2)
 \end{equation}
 Hence, at the lowest order, a non-zero branch of solution with $\rho \neq 0$ exists if and only if the term into the bracket vanishes, which gives rise to the instability region of this problem. Separating real and imaginary parts of this term, and eliminating the angle $\phi$ thanks to the relationship $\cos^2 2 \phi + \sin^2 2 \phi = 1$, yields:
 \begin{equation}\label{eq:branchpoint}
     \Omega_{a,b} = 2\omega_1\sqrt{1-\xi_1^2} \pm 2 \sqrt{f_2^2 - \xi_1^2 \omega_1^2}.
 \end{equation}
 This equation gives the two limit points $\Omega_{a,b}$ where the two branches of periodic solutions are connected to the solution at rest $\rho=0$, {\it i.e.} the two period-doubling bifurcation points. 

 The analysis can be even pushed further by considering the term $f_3 z_1 z_3 z_4 $ in addition. Since $z_3=\e^{\iu\omega t}$ and $z_4=\e^{-\iu\omega t}$, this trivially resonant monomial simply reads $f_3 z_1$ and will thus have a direct consequence on the linear term and on the oscillation frequency, which is perturbed. Considering an order 2 on the non-autonomous coordinate will thus change the prediction of the location of the bifurcation points given in Eq.~\eqref{eq:branchpoint}, which needs to be rewritten by changing $\omega_1$ to $\omega_1+f_3$. All those predictions will be assessed thanks to the numerical results, underlining the predictive capability and the wealth of the ROM.

Figs.~\ref{fig:param1}-\ref{fig:param2} show the results obtained with the ROM in the parametric excitation case, and investigates the convergence of the ROM by increasing the order of truncation. Differently to the previous cases that have been shown (primary and superharmonic response), a Rayleigh damping matrix $\alpha \bfM + \beta \bfK$ has here been used with $\alpha=-3.81\,10^{-2}s^{-1}$ and $\beta=2\,10^{-1}s$. This change is motivated by the fact that, for this parametric resonance, the nonlinear damping behaviour brought about by the slave modes is too small if one considers mass-proportional damping, leading to bifurcated branches of solutions that do not close and go to very large amplitude values. In this specific case, the nonlinear damping term, which is created from the linear damping values of the slave modes~\cite{TOUZE:JSV:2006,artDNF2020}, plays a major role in producing bifurcated branches of solutions with non-diverging values. In order to increase the linear damping values of the slave modes whilst keeping the linear damping of the master mode the same, we used a combination of stiffness-proportional and mass-proportional damping.

Fig.~\ref{fig:param1} shows the bifurcated solution branches of the parametric resonance for increasing odd orders (3,5,7 and 9) of the autonomous expansion, while keeping the non-autonomous truncation to 1. Interestingly, the y-axis represents the amplitude of the physical displacement of the center point of the beam $u_m$. Note that this value is small but not vanishing. Indeed, since $\omega \simeq 2\omega_1$, there is a residual non-resonant linear excitation, far from resonance, which makes the displacement solution non-zero here. Whereas the ROM coordinates $z$ exactly bifurcates from $\rho=0$, it is the nonlinear mapping that relates physical coordinates to normal coordinates, which allows retrieving this non-zero (though small) value of the displacement. The bifurcated branches of solutions then show the classical portrayal of a parametrically excited system, with a hardening behaviour obtained with the parameters selected here. Again, stiffness-proportional damping was needed here to reach a maximal amplitude in the bifurcated branches. Simulations with mass-proportional damping (not reported here for the sake of brevity) had branches that do not close and then diverge to infinity, underlining the importance of the damping of the slave modes which is accurately reproduced by the ROM and offering accurate and quantitative predictions.

Finally, Fig.~\ref{fig:param2} shows the convergence of the ROMs when increasing the order of the non-autonomous terms. As awaited from the analysis, an important shift of the solution is obtained when increasing the order of the non-autonomous coordinates to 2, since including the term $f_3 z_1 z_3 z_4$ in the analysis. As announced, taking this term into account leads to a slight frequency shift such that the location of the bifurcation points is slightly modified. Higher order terms then show minimal effects and the convergence is finally reached with order $\mathcal{O}(9,3)$.

\subsection{Computational costs}
All the computations presented in the preceding sections were performed utilizing an Intel i9-12900K processor with 32 GB of RAM.
The mesh employed for the parallel plate formulations, both for the reduced model and the full one, counts 5834 degrees of freedom.
With this mesh, the {\sc ppf-fom} model with the Newmark scheme requires approximately 1 hour per computed point.
The computational time needed to construct the ROM models depends on the selected order of expansion. In the case of the {\sc mppf-rom} formulation, the following computation times were recorded:
\begin{itemize}
    \item Order $\mathcal{O}(7,1)$: less than 1 minute;
    \item Order $\mathcal{O}(9,1)$: 3 minutes;
    \item Order $\mathcal{O}(9,5)$: 25 minutes.
\end{itemize}
Once the reduced model has been retrieved, the continuation procedure for each curve takes from 10 seconds to 1 minute.
The specific timing within this range depends on the number of points computed on the curve.

\section{Conclusions}
\label{app:conclusions}

In this paper, the direct parametrisation method for invariant manifolds (DPIM) has been applied to a fully coupled multiphysics problem involving the nonlinear vibrations of deformable structures under the influence of an electrostatic field, hence significantly extending the realm of applications that can be dealt with this reduction technique.

The main difficulty was related to the expression of the nonlinear forces of the electromechanical problem, which are nonpolynomial, thus asking for a special treatment in order to be then automatically treated in the arbitrary order expansions issued from the DPIM.  To overcome this challenge, a new mixed fully Lagrangian formulation has been proposed, that exclusively incorporates explicit polynomial nonlinearities. This was achieved by rewriting the coupled problem in its original configuration, and by introducing an auxiliary field strictly related to the electric field.



The governing equations have been semi-discretised using the standard finite element procedures, resulting in a system of differential-algebraic equations compatible with the general derivation of the DPIM for non-autonomous problems discussed in~\cite{vizza2023superharm}.

The proposed formulation is fully general and applicable to arbitrary geometries.
For the purpose of validation, we have restrained ourselves to a 2D implementation of the standard, but challenging benchmark of a clamped-clamped beam facing a parallel electrode,
where well-known and accurate simplified expressions for the electrostatic forces can be used for comparison.

This simplification enabled validation of the reduced-order models against full-order simulations obtained through direct time integration. Numerical results have shown a remarkable performance both in terms of accuracy and wealth of nonlinear effects that can be accounted for, including primary resonances, with the transition of the frequency response curves from hardening to softening, as well as secondary resonances. The convergence of the ROM has been also addressed and its excellent predictive capacity is used to give insights into the nonlinear dynamical phenomena observed.

Application of this nonlinear reduction technique paves the way for the generation of efficient and reliable 
ROMS with excellent accuracy and efficient predictive capacities for electrostatically actuated resonating Micro Electro Mechanical Systems,
and provides a much-needed tool for the design and optimization of a whole new family of devices with unprecedented performances.
This research contributes significantly to advancing the understanding and application of DPIM in tackling complex multiphysics problems, particularly those arising in the domain of MEMS.


\section*{Acknowledgments}
Attilio Frangi and Alessio Colombo acknowledge the support of ST Microelectronics under the framework of the Joint Research Platform with the Politecnico di Milano, STEAM.

\bibliographystyle{unsrt}
\bibliography{biblioROM}


\section*{Appendix}
\begin{appendix}

\section{Virtual power of internal stresses}
\label{app:MMCdetail}

In this section, we briefly summarize the steps left implicit in 
Section~\ref{eq:PPVmech} that allows reformulating the terms associated to the power of internal stresses in Eq.~\eqref{eq:PPVcurrent}.
The Cauchy stress tensor $\bfsig$ is replaced by the
second Piola-Kirchhoff stress $\bfsig\pk$ defined as
\[
\bfsig\pk = J \bff^{-1}\cdot\bfsig\cdot\bff^{-T}
\]
where $\bff[\bfu]=\bfone+\nabla\bfu$ is the transformation gradient
and $J$ is the jacobian of the transformation, so that
\[
\int_{\Omega} \bfsig:\grad\tilbfu\,\dOmega = 
\int_{\Omega_0} (\bff[\bfu]\cdot\bfsig\pk[\bfu]):\nablaT\tilbfu\,\dOmega_0
\]
%
Applying standard rules of tensor calculus:
\begin{align}
\label{eq:Pi}
\int_{\Omega_0} (\bff[\bfu]\cdot\bfsig\pk[\bfu]):\nablaT\tilbfu\,\dOmega_0
& =
\int_{\Omega_0} \bfsig\pk[\bfu]:\tsym(\nablaT\tilbfu\cdot\bff)\,\dOmega_0
\nnum
& =
\int_{\Omega_0} \bfsig\pk[\bfu]:(\bfeps[\tilbfu] + \bfens[\bfu,\tilbfu])\,\dOmega_0
\end{align}
where:
\begin{align}
\tsym(\nablaT\tilbfu\cdot\bff) & = \frac{1}{2}\left(\nabla\tilbfu + \nablaT\tilbfu + \nablaT\tilbfu\cdot\nabla\bfu + \nablaT\bfu\cdot\nabla\tilbfu \right)
\end{align}
and we have defined the operators
\begin{align}
\label{eq:eps}
\bfeps[\bfa] & =\frac{1}{2}(\nabla\bfa+\nablaT\bfa) 
\\
\label{eq:ens}
\bfens[\bfa,\bfb] & = \frac{1}{2}(\nablaT\bfa\cdot\nabla\bfb+\nablaT\bfb\cdot\nabla\bfa). 
\end{align}
Note also that: 
\[
\bfe[\bfu]=\bfeps[\bfu] + \frac{1}{2}\bfens[\bfu,\bfu]
\]
so that:
\begin{align}
\label{eq:split}
\bfsig\pk[\bfu]:(\bfeps[\tilbfu] + \bfens[\bfu,\tilbfu]) = & \bfeps[\tilbfu]:\cA:\bfeps[\bfu]+
\\
& \frac{1}{2} \bfeps[\tilbfu]:\cA:\bfens[\bfu,\bfu] +   \bfens[\tilbfu,\bfu]:\cA:\bfeps[\bfu] + 
\nnum
& \frac{1}{2} \bfens[\tilbfu,\bfu]:\cA:\bfens[\bfu,\bfu]
\nonumber
\end{align}
The first term is linear, the second and the third are quadratic and the last term is cubic. 

In many applications, the solution is sought as $\bfu_0+\bfu$
where $\bfu_0$ is known and $\bfu$ is the unknown increment.
Indeed, this is a slight abuse of notation, since now $\bfu$ denotes the increment and not the total
displacement, but simplicity will benefit from this choice.
It is then important to identify the new contributions in the integrand of Eq.\eqref{eq:Pi}:
\[
\bfsig\pk[\bfu]:(\bfeps[\tilbfu] + \bfens[\bfu,\tilbfu]) 
\]
The factors can be expanded as:
\begin{align}
\bfe[\bfu_0+\bfu] & =\bfe[\bfu_0]+\bfeps[\bfu]+\bfens[\bfu_0,\bfu]+
\frac{1}{2}\bfens[\bfu,\bfu]
\\
\bfeps[\tilbfu] + \bfens[\bfu_0+\bfu,\tilbfu] & =
\bfeps[\tilbfu] + \bfens[\bfu_0,\tilbfu]+\bfens[\bfu,\tilbfu]
\end{align}
so that the product yields a constant term:
\[
\bfsig\pk[\bfu_0]:(\bfeps[\tilbfu] + \bfens[\bfu_0,\tilbfu]) 
\]
a linear term:
\[
\bfsig\pk[\bfu_0]:\bfens[\bfu,\tilbfu]+
(\bfeps[\bfu] + \bfens[\bfu_0,\bfu]):\cA:(\bfeps[\tilbfu] + \bfens[\bfu_0,\tilbfu]) 
\]
a quadratic term:
\[
(\bfeps[\bfu]+\bfens[\bfu_0,\bfu]):\cA:\bfens[\bfu,\tilbfu]+
\frac{1}{2}\bfens[\bfu,\bfu]:\cA:(\bfeps[\tilbfu] + \bfens[\bfu_0,\tilbfu])
\]
and finally a cubic one:
\[
\frac{1}{2}\bfens[\bfu,\bfu]:\cA:\bfens[\bfu,\tilbfu]
\]
The FEM discretisation of the internal power leads to:
\begin{align}
\label{eq:dpimmech}
\tilde{\bfU}^\tT\left(\bfF+\bfK\bfU+\bsnG(\bfU,\bfU)+
\bsnH(\bfU,\bfU,\bfU)\right) 
\end{align}
where the matrices can be computed directly from the decomposition above.

\section{Discretisation of the {\sc mppf} formulation}
\label{app:mafdiscr}

Let us consider the simplified approach for parallel plates of Section~\ref{sec:PPF}
formulated by the coupled continuous equations \eqref{eq:MAFmech}
and \eqref{eq:MAFcomp}.
Following the same procedure as for the {\sc mgf} formulation,
let us denote with $\bfu_0,\psi_0$ the static solution 
of the {\sc mppf}  eqs.\eqref{eq:MAFmech}-\eqref{eq:MAFcomp}
when only $V_\tDC$ is acting
and with $\bfu_0+\bfu,\psi_0+\psi$ the total solution under the combined action
of $V_\tDC$ and $V_\tAC$.
Ignoring the constant terms which enter only in the 
computation of the static solution
the virtual power of electrostatic forces in 
Eq.\eqref{eq:MAFmech} becomes:
\begin{align}
\label{eq:MAFmechD}
\cP_e=
\underbrace{
\veps_0\int_{\tU_0} \psi_0\psi\tilu_2\,\dS_0
}_{\text{linear}}
+
\underbrace{
\frac{\veps_0}{2}\int_{\tU_0} \psi^2\tilu_2\,\dS_0
}_{\text{quadratic}}
\end{align}
so that the semi-discretised equation for the mechanical equilibrium reads,
also using Eq.\eqref{eq:dpimmech}:
\begin{align}
\label{eq:MAF1}
\bfM\ddot{\bfU}+\bfK\bfU+\bsnG(\bfU,\bfU)+\bsnH(\bfU,\bfU,\bfU)
 = \bfR_\psi\bfPsi+\bsnR_{\psi\psi}(\bfPsi,\bfPsi)
\end{align}
The compatibility equation \eqref{eq:MAFcomp} has a right hand side which generates the forcing while the left hand generates the terms:
\begin{align}
\label{eq:MAFcompD}
\underbrace{
\int_{\tU_0}(\psi(g-u_{02})-\psi_0 u_2)\tilpsi\,\dS_0
}_{\text{linear}}
+
\underbrace{
\int_{\tU_0}-\psi u_2\tilpsi\,\dS_0
}_{\text{quadratic}}
\end{align}
The final form of the discretised compatibility equation is:
\begin{align}
\label{eq:MAF2}
\bfC_u\bfU+\bfC_{\psi}\bfPsi+ \bsnC_{\psi u}(\bfPsi,\bfU)
= \bfF_{\psi}V_{\tAC}\sin(\omega t)
\end{align}
where the right-hand side vector comes from:
\[
\tilde{\bfPsi}^\tT \bfF_{\psi}
=
\int_{\tU_0} \tilpsi\,\dS
\]
The system of Eqs.~\eqref{eq:MAF1}-\eqref{eq:MAF2}
is a Differential Algebraic System which can be reduced via the DPIM as 
described in Section \ref{sec:dpim} for the equivalent DAE
stemming from the {\sc mgf} of Section \ref{sec:GF}.

\section{Elements of differential geometry of a surface}
\label{app:surface}

Let us consider a surface $S_0$
embedded in a three-dimensional Euclidean space $x_1,x_2,x_3$ 
and mapped onto the
parameter space $a_1,a_2$ (e.g. the parameters of a FEM discretisation).

We now introduce $\bfg_{\alpha}$, the two covariant base vectors:
\begin{align}
\label{a1}
\bfg_{\alpha}=\bfx_{,\alpha} \, (= \der{\bfx}{a_{\alpha}})
\end{align} 
which are in general non-orthogonal, non-unitary vectors tangent to the surface.
The vector product of the two covariant vectors yields the normal:
\begin{align}
\label{a2}
\bfg_1\wedge\bfg_2=  \bfn J_S
\quad \text{where} \quad
J_S=\|\bfg_1\wedge\bfg_2\|
\end{align} 
The differential surface element can be expressed as:
\begin{align}
\label{eq:ds0}
\dS_0 =  J_S\,\da_1\da_2 
\end{align} 
and in particular, one has:
\begin{align}
\label{eq:Nds0}
\bfn\dS_{0}=\bfg_1\wedge\bfg_2\,\da_1\da_2 
\end{align}
%

Let us now suppose that the surface deforms due to displacements $\bfu$ with
$\bfy=\bfx+\bfu$ being the actual coordinates. 
Let $\hat{\bfg}_{\alpha}$ denote the ``actual'' covariant vectors
\begin{align}
\label{actual_covariant_vectors}
\hat{\bfg}_{\alpha}=\bfy_{,\alpha}=
\bfg_{\alpha}+\bfu_{,\alpha} \, 
(=\bfg_{\alpha}+\nabla\bfu\dod\bfg_{\alpha})
\end{align} 
The vector product of the new covariant vectors gives the actual normal:
\begin{align}
\label{a2_actual}
\hat{\bfg}_1\wedge\hat{\bfg}_2 = \hat{J}_S \hat{\bfn}
\end{align} 
We first compute the expansion of the vector product $\hat{\bfg}_{1}\wedge\hat{\bfg}_{2}$:
\begin{align}
\label{eq:g1g2}
\hat{\bfg}_1\wedge\hat{\bfg}_2 & = (\bfg_1+\bfu_{,1})\wedge(\bfg_2+\bfu_{,2}) 
\\
& = \bfg_1\wedge\bfg_2+ \bfg_1\wedge\bfu_{,2}+\bfu_{,1}\wedge\bfg_2 + 
\bfu_{,1}\wedge\bfu_{,2}
\nonumber
\end{align} 
showing that the surface element $\hat{J}_S \hat{\bfn}$ 
admits an exact quadratic expansion in $\bfu$:
\begin{align}
\label{eq:njs}
\hat{J}_S \hat{\bfn}  = J_S \bfn +  \bfg_1\wedge\bfu_{,2}+\bfu_{,1}\wedge\bfg_2 + 
\bfu_{,1}\wedge\bfu_{,2}
\end{align} 
In particular, one has:
\begin{align}
\label{eq:Nds}
\hat{\bfn}\dS=\hat{\bfg}_1\wedge\hat{\bfg}_2\,\da_1\da_2 =
&
(J_S \bfn +  \bfg_1\wedge\bfu_{,2}+\bfu_{,1}\wedge\bfg_2 + 
\bfu_{,1}\wedge\bfu_{,2})\da_1\da_2
\nnum
=
& (\bfn +  \frac{1}{J_S}\bfg_1\wedge\bfu_{,2}+\bfu_{,1}\wedge\bfg_2 + 
\bfu_{,1}\wedge\bfu_{,2})\dS_0
\end{align}

\end{appendix}

\end{document}